\documentclass[11pt]{article}
\usepackage{pdflscape}
\usepackage{latexsym,amsmath,amssymb,mathrsfs,graphicx,hyperref,longtable,mathpazo,soul,color,amsthm,mathtools,xcolor,multirow}
\hypersetup{pdftex,colorlinks=true,linkcolor=blue,citecolor=blue,filecolor=blue,urlcolor=blue}

\usepackage{url}
\usepackage{caption}
\def\tr{{\raise0pt\hbox{$\scriptscriptstyle\top$}}}

\usepackage{longtable}

\usepackage{rotating}
\newtheorem{theorem}{Theorem}[section]
\newtheorem{corollary}[theorem]{Corollary}

\newtheorem{proposition}[theorem]{Proposition}

\newtheorem{problem}[theorem]{Problem}
\newtheorem{definition}[theorem]{Definition}
\newtheorem{question}[theorem]{Open Question}

\newtheorem{conjecture}[theorem]{Conjecture}

\numberwithin{equation}{section}
\numberwithin{table}{section}
\title{Generalised Fermat equation: a survey of solved cases}
\author{Ashleigh Ratcliffe and Bogdan Grechuk}
\begin{document}
	\maketitle
	
\begin{abstract}
	Generalised Fermat equation (GFE) is the equation of the form $ax^p+by^q=cz^r$, where $a,b,c,p,q,r$ are positive integers. If $1/p+1/q+1/r<1$, GFE is known to have at most finitely many primitive integer solutions $(x,y,z)$. A large body of the literature is devoted to finding such solutions explicitly for various six-tuples $(a,b,c,p,q,r)$, as well as for infinite families of such six-tuples. This paper surveys the families of parameters for which GFE has been solved. Although the proofs are not discussed here, collecting these references in one place will make it easier for the readers to find the relevant proof techniques in the original papers. 
	Also, this survey will help the readers to avoid duplicate work by solving the already solved cases.
\end{abstract}

\tableofcontents

\section{Introduction}

\subsection{Generalized Fermat equation}

Around 1637, Pierre de Fermat wrote in the margin of his private copy of Arithmetica that for any integer $n\geq 3$ the equation
\begin{equation*}
x^n + y^n = z^n
\end{equation*}
has no solutions in positive integers. This statement became known as Fermat's Last Theorem (FLT), and
has been one of the most famous open problems in mathematics for over $350$ years. After
cancelling the common factors if needed, it suffices to prove the non-existence of positive integer
solutions such that $\gcd(x,y,z)=1$. The FLT has been finally proved by Andrew Wiles\footnote{With part of
the proof delegated to a joint paper with Richard Taylor~\cite{MR1333036}.}~\cite{wiles1995modular} in
1995, and this proof is widely recognized as one of the most monumental achievements of modern
mathematics.

Even before FLT had been proven, many researchers started to investigate special cases of the more general equation
\begin{equation}
\label{eq:genfermat}
 ax^p+by^q=cz^r,
\end{equation} 
where $x,y,z$ are integer variables, while positive integers $p,q,r$ and integers $a,b,c$
are parameters. If $abc=0$, then \eqref{eq:genfermat} has at most two monomials and is easy to
solve, see~\cite[Section 1.6]{mainbook} and~\cite[Section 1.16]{wilcox2024systematic}, so from now on
we will assume that $abc\neq 0$. If \eqref{eq:genfermat} is solvable in integers $(x,y,z)\neq (0,0,0)$, then,
under a minor condition\footnote{Specifically, under the condition that at least one of the integers
$p,q,r$ is coprime with another two.} on $(p,q,r)$, it has infinitely many integer solutions,
and all these solutions are easy to describe, see~\cite[Section 3.4.2]{MR4620765}
and~\cite[Proposition 4.53]{mainbook}. However, the problem of describing \emph{primitive} integer
solutions to \eqref{eq:genfermat}, that is, ones satisfying $\gcd(x,y,z)=1$, turns out to be much more
interesting and delicate. 

The triple $(p,q,r)$ is called the signature of Eq.~\eqref{eq:genfermat}. If $\min\{p,q,r\} = 1$ then
\eqref{eq:genfermat} is trivial, see~\cite[Section 1.5.1]{mainbook} and~\cite[Section
1.12]{wilcox2024systematic}, so we may assume that $\min\{p,q,r\}\geq 2$. The theory of primitive solutions to
\eqref{eq:genfermat} crucially depends on whether the quantity $1/p  +  1/q   +  1/r$ is greater than
$1$, equal to $1$, or less than $1$. In the first case, we have the following
result proved by Beukers~\cite{MR1487980} in 1998.

\begin{theorem}
\label{eq:Beukers1998}
 Assume that $1/p+1/q+1/r>1$ and $\min\{p,q,r\}\geq 2$. Then if \eqref{eq:genfermat} has at least one primitive solution with $xyz\neq 0$, then it has infinitely many such solutions. In this case, there exists a finite number of parametrized solutions of \eqref{eq:genfermat} with two parameters such that all primitive solutions can be obtained by specialization of the parameters to integer values. Moreover, there is an algorithm that computes these parametrized solutions.
\end{theorem} 

Theorem \ref{eq:Beukers1998} covers the signatures $(2,2,k)$ for $k\geq 2$, $(2,3,3)$, $(2,3,4)$,
$(2,3,5)$, and their permutations. As a nice example of its application,
see~\cite{edwards2005platonic} for the complete description of primitive integer solutions to the
equation $x^2+ y^3= z^5$. 

In the case $1/p+1/q+1/r=1$, which consists of permutations of $(2,3,6)$, $(2,4,4)$ and $(3,3,3)$,
Eq.~\eqref{eq:genfermat} reduces to the problem of computing rational points on certain elliptic
curves. This problem is open in general, but there are algorithms that work well for all specific
equations with not-too-large coefficients $a,b,c$. For example, as early as in 1954,
Selmer~\cite{MR41871,MR67131} solved \eqref{eq:genfermat} with $(p,q,r)=(3,3,3)$ for all $a,b,c$
satisfying $1\leq | abc|  \leq 500$. With modern computers and improved algorithms one may proceed much further if
needed.  

The most interesting and challenging case of Eq.~\eqref{eq:genfermat} is $1/p  +  1/q   +  1/r   <   1$. In this case, Darmon
and Granville~\cite{DOI:10.1112/blms/27.6.513} proved the following 
fundamental result.

\begin{theorem}
\label{th:genfermat}
 For any given positive integers $p,q,r,a,b,c$ satisfying $1/p  +  1/q   +  1/r   <   1$,  the generalized Fermat Eq.~\eqref{eq:genfermat} has only finitely many primitive integer solutions.
\end{theorem}

\begin{table}
\begin{center}
  \caption{\label{tb:pqrsigs}Triples $(p,q,r)$ for which equation \eqref{eq:fermatpqr} has been solved, up to permutations (triples marked with * do not include permutations).}
\begin{tabular}{|c|c|}
\hline
   Triple & Solved  \\ \hline \hline
   $(2,3,n), 7 \leq n \leq 12, $\footnotemark $n=15$ & \cite{2007ToXa,BRUIN_1999,+2003+27+49,BRUIN2005179,dahmen2008,brown2010primitive,MR3095226,MR4036449,SiksekStoll14}  \\\hline
   $(2, 4, n), n \geq 5$ & \cite{bennett2010diophantine,MR2075481,BRUIN_1999,Bennett_Skinner_2004,+2003+27+49}  \\\hline
   $(2,6,n), n \geq 3$ & \cite{bennett2012multi,BRUIN_1999,miscellany} \\\hline
   $(2,2n,3)^*, 3 \leq n \leq 10^7$ & \cite{chen2008equation,dahmen2011refined,dahmen2008,MR3095226} \\\hline 
   $(2,2n,m)^*, m \in\{6,9,10,15,21\}, n\geq 2$ & \cite{miscellany,MR4418449}  \\\hline
           $(4,2n,3)^*, n \geq 2$ & \cite{miscellany} \\\hline
       $(3,4,5)$ & \cite{siksek2012partial}  \\\hline
   $(3,3,n), 3 \leq n \leq 10^9$ & \cite{chen2009perfect,kraus1998equation,10.1007/10722028_9,dahmen2008,MR3526934} \\\hline
   $(3,3,2n), n \geq 2$ & \cite{miscellany}  \\\hline
   $(3,6,n), n \geq 3$ & \cite{miscellany} \\\hline
   $(5,5,n), n \in \{7,19\}$ & \cite{dahmen2014perfectpowersexpressiblesums} \\\hline
   $(7,7,5)$ & \cite{dahmen2014perfectpowersexpressiblesums} \\\hline 
   $(n, n, 2), n \geq 5$ & \cite{Darmon1997,poonen1998some} \\\hline
   $(n, n, 3), n \geq 3$ & \cite{Darmon1997,poonen1998some}\footnotemark  \\\hline
   $(n, n, n), n \geq 3$ & \cite{wiles1995modular,MR1333036} \\\hline
   $(2j,2k,n)^*, j,k\geq 5$ prime, $n \in \{3,5,7,11,13\}$ & \cite{anni2016modular,10.5565/PUBLMAT6722309}  \\\hline
   $(2j,2k,17)^*, j,k \neq 5$ primes & \cite{MR4418449} \\\hline
   $(2n, 2n, 5)^*, n \geq 2$ & \cite{JTNB_2006} \\\hline
   $(2n,2n,17)^*$ & \cite{MR4418449} \\\hline
   $(3j, 3k, n)^*, j, k \geq 2, n \geq 3$ & \cite{kraus1998equation} \\
\hline
\end{tabular}
\end{center}
\end{table}
\footnotetext[3]{The proof for $n=11$ in \cite{MR4036449} is conditional on the generalized Riemann hypothesis. The case $n=12$ follows from the analysis of case $(2,3,6)$, which reduces to a rank $0$ elliptic curve.}
\footnotetext[4]{It was noted on Mathoverflow (\url{https://mathoverflow.net/questions/488724/}) that \cite{Darmon1997} does not provide details for the case $(3,4,4)$, this case is proven in \cite[Proposition 14.6.6]{MR2312338}.}

The proof of Theorem \ref{th:genfermat} is ineffective and does not provide an algorithm for actually listing the primitive solutions to \eqref{eq:genfermat} for a given $p,q,r,a,b,c$. This problem is the topic of current intensive research. The case $a=b=c=1$, that is, equation
\begin{equation}
\label{eq:fermatpqr}
 x^p+y^q=z^r,
\end{equation}
is known as the Fermat--Catalan equation, and has been particularly well-studied. A number of survey
papers exist which are devoted to this equation, see, for example,~\cite{miscellany,Kraus1999},
Section 1.2 in~\cite{MR4122899}, and  wikipedia page
\url{https://en.wikipedia.org/wiki/Beal\_conjecture}. We provide a summary\footnote{Recently,
Bartolom\'e and Mih{\u a}ilescu~\cite{bartolome2021semilocal} claimed to solve \eqref{eq:fermatpqr} for
all signatures $(p,p,r)$ such that $p,r$ are primes, $p\geq 5$, and $r>3\sqrt{p\log_2 p}$. To the best
of our knowledge, this paper has not been peer reviewed yet.} of solved signatures $(p,q,r)$ in 
Table \ref{tb:pqrsigs}. There is also some recent work with partial results for other signatures, see
e.g~\cite{MR4205757,MR3217670,miscellany,MR4306226,MR3690598,signature2019multi,10.1007/10722028_9,MR2652585,chen2022modularapproachfermatequations,chen2009perfect,dahmen2008,dahmen2011refined,dahmen2014perfectpowersexpressiblesums,dahmen2024generalized,dieulefait2005modular,MR3493372,MR4036449,kraus1998equation,madriaga2024hypergeometricmotivesgeneralizedfermat,MR4648747,ThesisPutz,MR4580454}.
  Famous Fermat--Catalan conjecture predicts that \eqref{eq:fermatpqr} has only finitely many solutions in positive integers $(p,q,r,x,y,z)$ such that $1/p+1/q+1/r<1$ and $\text{gcd}(x,y,z)=1$. In fact, the only known such solutions are listed in 
Table \ref{tb:knownsol}, and failing to find other solutions despite extensive
efforts~\cite{norvig2015beal,sikora2024fermat} suggests that this list might be complete.  One may
note that all exponents $(p,q,r)$ in 
Table \ref{tb:knownsol} have $\min(p, q, r) = 2$. In 1997, Beal~\cite{MR1488570} offered a million-dollar prize
for the proof or disproof that \eqref{eq:fermatpqr} has no solutions in coprime positive integers
$(x,y,z)$ when $\min\{p,q,r\}\geq 3$.

\begin{table}
\begin{center}
  \caption{\label{tb:knownsol}Known primitive positive integer solutions to~\eqref{eq:fermatpqr} with $1/p+1/q+1/r<1$, up to exchange of $(x,p)$ with $(y,q)$.}
\begin{tabular}{|c|c|c|c|}
\hline

   $(p,q,r)$ & $(x,y,z)$ & $(p,q,r)$ & $(x,y,z)$  \\\hline \hline 
   $(p,3,2)$ & $(1,2,3)$ &  $(8,2,3)$ & $(33,1549034,15613)$  \\\hline
   $(5,2,4)$ & $(2,7,3)$ & $(3,2,7)$ & $(1414,2213459,65)$  \\\hline
   $(3,2,9)$ & $(7,13,2)$ & $(3,2,7)$ &  $(9262,15312283, 113)$ \\\hline
   $(7,3,2)$ & $(2,17,71)$ & $(7,3,2)$ & $(17,76271,21063928)$ \\\hline
   $(5,4,2)$ & $(3,11,122)$ & $(8,3,2)$ & $(43, 96222,30042907)$ 
\\\hline
\end{tabular}
\end{center}
\end{table}

Comparison of 
Table \ref{tb:pqrsigs} with 
Table \ref{tb:triples} below shows that the only cases of Eq.~\eqref{eq:fermatpqr} with $(p,q,r)$ not-all-primes that remain to be considered are $(p,q,r)=(2,5,9)$ and $(p,q,r)=(2,3,25)$.  In the case when all $(p,q,r)$ are primes, 
Table \ref{tb:pqrsigs} covers the triples $(2,3,7)$, $(5,5,7)$, $(5,5,19)$, $(7,7,5)$, $(3,3,p)$ for $p\leq 10^9$, $(p,p,2)$, $(p,p,3)$ and $(p,p,p)$ for all $p$, as well as $(2,3,11)$ conditional on the generalized Riemann hypothesis. Hence, the triple $(p,q,r)$ with the smallest sum $p+q+r$ for which Eq.~\eqref{eq:fermatpqr} is currently open is $(p,q,r)=(2,5,7)$. The smallest triple for which Beal's conjecture remains open is $(p,q,r)=(3,5,7)$. 

Recently, Sikora~\cite{sikora2024fermat} reported that Eq.~\eqref{eq:fermatpqr} has no solutions in
coprime positive integers other than ones listed in 
Table \ref{tb:knownsol} in the range $z^r<2^{71}$. Restricting the search to only exponents not covered in 
Table \ref{tb:pqrsigs} allowed us to significantly speed up the search and reach the bound $2^{100}$.

\begin{proposition}
\label{prop:onlysol}
 The only solutions to \eqref{eq:fermatpqr} in coprime positive integers such that $z^r\leq 2^{100}$ are those listed in 
Table \ref{tb:knownsol}.
\end{proposition} 

While there are good surveys of Eq.~\eqref{eq:fermatpqr}, a similar survey for the more general Eq.~\eqref{eq:genfermat} is not available. There are many dozens of papers studying this equation, and
each of them cites just a small portion of others. As a result, some equations of the form
\eqref{eq:genfermat} have been solved multiple times. For example, equation $x^4+3y^4=z^3$ has been
solved in 1994 by Terai~\cite{MR1288426}, and then again in 2019 by
S\"oderlund~\cite{soderlund2019some}.

The aim of this survey is to summarize the solved cases of the generalized Fermat Eq.~\eqref{eq:genfermat} with $| abc|  > 1$. This will help the readers to avoid duplicate work by solving the already solved cases. To keep the survey relatively short, we will not discuss proofs. Despite this, we hope that collecting all these references in one place should make it much easier for the readers to find the proof techniques in the original papers. 

In the rest of the introduction, we discuss some easy reductions, elementary special cases, and state some conjectures. Sections \ref{sec:special} to \ref{sec:pqr} discuss for which parameters $(a,b,c,p,q,r)$ Eq.~\eqref{eq:genfermat} has been solved. The solved cases are then summarized in the final table (Table \ref{tbl10}) in Section \ref{sec:summary}.

\subsection{Easy reductions and elementary special cases}

We first observe that we may assume that the integers $a,b$ and $c$ in \eqref{eq:genfermat} are all positive. Indeed, if an equation $Ax^p+By^q+Cz^r=0$ has $A,B,C$ all of the same sign, then it has no solutions other than $(x,y,z)=(0,0,0)$ if all $p,q,r$ are even, and otherwise we can change the sign of one of $A,B,C$ by a transformation of one of the variables to their negative. Hence, we may assume that not all $A,B,C$ are of the same sign, and then the equation can be rearranged to the form \eqref{eq:genfermat} with positive $a,b,c$. 

We are interested in finding primitive solutions to \eqref{eq:genfermat}, that is, ones satisfying $\gcd(x,y,z)=1$. Because solutions to \eqref{eq:genfermat} with $xyz=0$ are easy to describe, it suffices to find primitive solutions in non-zero integers. We will call such solutions \emph{non-trivial}. Further, if a prime $s$ is a common factor of, say, $x$ and $y$, then $\gcd(x,y,z)=1$ implies that $\gcd(s,z)=1$, but then \eqref{eq:genfermat} implies that $s^{\min\{p,q\}}$ is a divisor of $c$. Continuing this way, we obtain the following observation.

\begin{proposition}
\label{prop:paircop}
 Assume that $c$ has no divisors of the form $s^{\min\{p,q\}}$ for prime $s$, $b$ has no divisors of the form $s^{\min\{p,r\}}$ for prime $s$, while $a$ has no divisors of the form $s^{\min\{q,r\}}$ for prime $s$. Then $(x,y,z)$ is a primitive solution to \eqref{eq:genfermat} if and only if $x,y,z$ are pairwise coprime.
\end{proposition}

From now on, for all equations satisfying the conditions of Proposition \ref{prop:paircop}, we will use the terms ``primitive solution'' and ``solution in pairwise coprime integers'' interchangeably.

We next observe that if one of the integers $(p,q,r)$, say $r$, is not prime, and $1<d<r$ is any divisor of $r$, then change of variables $Z=z^{r/d}$ reduces \eqref{eq:genfermat} to an equation of the same form with exponents $(p,q,d)$. If $1/p+1/q+1/d<1$, then this equation has finitely many primitive integer solutions by Theorem \ref{th:genfermat}. If we list all such solutions $(x,y,Z)$, we can then check for which of them $z=Z^{d/r}$ is an integer. This implies that it is sufficient to solve \eqref{eq:genfermat} only in the case when either $p,q,r$ are all primes, or $(p,q,r)$ is a \emph{special} triple defined below.

\begin{definition}
A triple $(p,q,r)$ of positive integers is called \emph{special} if (i) not all $p,q,r$ are primes, (ii) $1/p+1/q+1/r<1$, and (iii) if $P,Q,R$ are positive divisors of $p,q,r$, respectively, such that $(P,Q,R)\neq (p,q,r)$, then $1/P+1/Q+1/R \geq 1$.  
\end{definition}

Every special triple is a permutation of one of the triples listed in 
Table \ref{tb:triples}.

\begin{table}
\begin{center}
 \caption{\label{tb:triples}Special triples $(p,q,r)$, up to permutations.}
\begin{tabular}{|c|c|c|c|}

\hline
   $(2,3,8)$ & $(2,3,25)$   & $(2,5,6)$ & $(3,3,9)$  \\\hline
   $(2,3,9)$  &$(2,4,r), r\geq 5 \text{ prime}$ & $(2,5,9)$ & $(3,4,4)$   \\\hline
   $(2,3,10)$ & $(2,4,6)$ & $(2,6,6)$    & $(3,4,5)$  \\\hline
   $(2,3,12)$ & $(2,4,8)$  & $(3,3,4)$  & $(4,4,4)$  \\\hline
   $(2,3,15)$   & $(2,4,9)$ &  $(3,3,6)$  & 
\\\hline
\end{tabular}
\end{center}
\end{table}
Some equations of the form \eqref{eq:genfermat} have no non-trivial primitive solutions for elementary reasons. For example, modulo $9$ analysis shows that if $(x,y,z)$ is any integer solution to the equation  $x^4 + y^4 = 3 z^3$, then all $x,y,z$ must be divisible by $3$. Hence, this equation has no non-trivial primitive solutions. More generally, we have the following trivial observation.

\begin{proposition}
\label{prop:local}
 If there exists any prime $s$ and positive integer $m$ such that all solutions $(x,y,z)$ to \eqref{eq:genfermat} modulo $s^m$ must have all $x,y,z$ divisible by $s$, then \eqref{eq:genfermat} has no non-trivial primitive solutions.
\end{proposition}

For equations covered by Proposition \ref{prop:local}, we say that they have no non-trivial primitive solutions by local obstructions. All such equations will from now be excluded from this survey.

As mentioned above, we may assume that either $p,q,r$ are all primes, or $(p,q,r)$ is a permutation of one of the special triples listed in 
Table \ref{tb:triples}. In the second case, at least one of the integers $p,q,r$ is composite. If, for example, $r$ is composite, and $1<d<r$ is any divisor of $r$, then change of variables $Z=z^{r/d}$ reduces \eqref{eq:genfermat} to the equation
\begin{equation}
\label{eq:genfermatred}
 ax^p+by^q=cZ^d.
\end{equation} 
Because $(p,q,r)$ is a special triple, we must have $1/p+1/q+1/d\geq 1$. If $1/p+1/q+1/d>1$, then, by Theorem \ref{eq:Beukers1998}, either \eqref{eq:genfermatred} has no non-trivial primitive solutions, or all such solutions can be covered by a finite number of parametrizations with two parameters. In the first case, \eqref{eq:genfermat} has no non-trivial primitive solutions as well. In the second case, we obtain that
\begin{equation}
\label{eq:genfermateasy}
 z^{r/d}=Z=P_i(u,v)  \text{ for some } i\in\{1,2,\dots,k\},
\end{equation} 
where $P_1,\dots,P_k$ are some polynomials in two variables $u,v$ with integer coefficients. In some cases, Eqs.~\eqref{eq:genfermateasy} are easy to solve, which leads to the resolution of \eqref{eq:genfermat}.

If $1/p+1/q+1/d=1$, then \eqref{eq:genfermatred} reduces to computing rational points on an elliptic
curve. If this curve has rank $0$, then there is a finite number of such points, and they can
be explicitly listed, see e.g.~\cite[Section 3.4.4]{mainbook} and~\cite[Section
3.21]{wilcox2024systematic}. Then \eqref{eq:genfermatred} has a finite number of primitive solutions,
and it is easy to check for which of them (if any) $z=Z^{d/r}$ is an integer. 

For some of the special triples in 
Table \ref{tb:triples} more than one of the exponents $p,q,r$ is composite, and some exponents may have more than one non-trivial divisor $d$. Hence, Eq.~\eqref{eq:genfermat} with special $(p,q,r)$ may be reducible to several equations of the form \eqref{eq:genfermatred}. If any of the reduced equations has finitely many primitive solutions, and these solutions can be computed, this immediately solves \eqref{eq:genfermat}. For example, this solves equation $x^4+2y^4=z^4$, which is the case
$(a,b,c,p,q,r)=(1,2,1,4,4,4)$ of Eq.~\eqref{eq:genfermat}, because the corresponding equation $x^4+2y^4=Z^2$ reduces to finding rational points on rank $0$ elliptic curve $Y^2=X^4+2$, where $Y=\frac{Z}{y^2}$, $X=\frac{x}{y}$. From now on, all Eqs.~\eqref{eq:genfermat} solvable by this method will be excluded from this survey.

In summary, below we will discuss the generalized Fermat's Eqs.~\eqref{eq:genfermat}, not covered by Proposition \ref{prop:local}, in which $a,b,c$ are positive integers with $abc>1$, and exponents $(p,q,r)$ are either (i) prime numbers satisfying $1/p+1/q+1/r<1$, or (ii) one of the special triples listed in 
Table \ref{tb:triples} such that all the corresponding Eqs.~\eqref{eq:genfermatred} have infinitely many primitive solutions. For each such equation, we will consider the following problem.

\begin{problem}
\label{prob:main}
 For a given $a,b,c,p,q,r$, list all primitive integer solutions $(x,y,z)$ to \eqref{eq:genfermat}. 
\end{problem}

For some equations, we will also discuss the following easier problem.
\begin{problem}
\label{prob:existence}
 For a given $a,b,c,p,q,r$, determine whether Eq.~\eqref{eq:genfermat} has any non-trivial primitive integer solutions.
\end{problem}

Because the solutions to \eqref{eq:genfermat} with $xyz=0$ are easy to find, the negative answer to Problem \ref{prob:existence} solves Problem \ref{prob:main} for the corresponding equation as well. However, if \eqref{eq:genfermat} has some easy-to-find non-trivial primitive integer solution, this automatically solves Problem \ref{prob:existence} with a ``Yes'' answer, without resolving Problem \ref{prob:main} for this equation.

\subsection{Some conjectures and open problems}

Famous and well-believed abc conjecture of Masser and Oesterl\'e~\cite{oesterle1988nouvelles} predicts
that for every real number $\epsilon>0$ there exist only finitely many triples $(A,B,C)$ of coprime
positive integers such that
\begin{equation}
\label{eq:abc}
 A+B=C \quad \text{and} \quad C > \mathrm{rad}(ABC)^{1+\epsilon},
\end{equation} 
where $\mathrm{rad}(n)$ denotes the product of the distinct prime factors of $n$. Darmon and
Granville~\cite[Section 5.2]{DOI:10.1112/blms/27.6.513} observed that the abc conjecture
implies the Fermat--Catalan conjecture about the finiteness of primitive integer solutions to
\eqref{eq:fermatpqr}, and remarked that the same argument may be extended to show that it also
implies the following conjectures.

\begin{conjecture}
\label{conj:abcfixed}
 For any integers $a,b,c$ such that $abc\neq 0$, Eq.~\eqref{eq:genfermat} has only finitely many solutions\footnote{Solutions to \eqref{eq:genfermat} that have $\min\{x,y,z\}=1$ and differ only by the exponent(s) of $1$ are assumed to be the same, e.g.~family $(x,y,z,p,q,r)=(1,1,1,p,q,r)$ for $(a,b,c)=(1,1,2)$ is considered as one solution.} in integers $(x,y,z,p,q,r)$ such that $\gcd(x,y,z)=1$, $p>0$, $q>0$, $r>0$, $1/p+1/q+1/r<1$.
\end{conjecture}

\begin{conjecture}
\label{conj:primesfixed}
 For any finite set $S$ of primes, Eq.~\eqref{eq:genfermat} has only finitely many solutions in integers \\ $(a,b,c,x,y,z,p,q,r)$ such that $\gcd(ax,by,cz)$ $=1$, $p>0$, $q>0$, $r>0$, $1/p+1/q+1/r<1$, and all the prime factors of $abc$ belong to set $S$.
\end{conjecture}

\begin{table}
\begin{center}
  \caption{\label{tab:pqrprimesols}Known positive integer solutions $(x,y,z,p,q,r,a,b,c)$ to \eqref{eq:open}, up to exchange of $(a,p,x)$ and $(b,q,y)$.}
\begin{tabular}{|c|c|c|}
\hline
$(a,b,c)$ & $(p,q,r)$ & $(x,y,z)$  \\
\hline \hline
   $(1,1,2)$ & $(p,q,r)$ & $(1,1,1)$ \\\hline
   $(1,2,1)$ & $(2,q,3)$ & $(5,1,3)$ \\\hline
   $(2,1,1)$ & $(2,q,5)$ & $(11,1,3)$ \\\hline
   $(1,1,2)$ & $(3,2,8)$ & $(7,13,2)$ \\\hline
   $(2,1,1)$ & $(4,4,3)$ & $(5,3,11)$ \\\hline
   $(1,1,2)$ & $(7,3,2)$ & $(3,5,34)$ \\\hline
   $(1,1,2)$ & $(3,q,2)$ & $(23,1,78)$ \\\hline
   $(1,2,1)$ & $(5,8,2)$ & $(7,3,173)$ \\\hline
   $(1,2,1)$ & $(5,4,2)$ & $(7,9,173)$ \\\hline
   $(2,1,1)$ & $(9,3,2)$ & $(3,19,215)$ \\\hline
   $(1,1,2)$ & $(2,q,4)$ & $(239,1,13)$ \\\hline
   $(1,1,2)$ & $(3,4,3)$ & $(237,43,203)$ \\\hline
   $(1,1,2)$ & $(3,13,2)$ & $(239,3,2761)$ \\\hline
   $(2,1,1)$ & $(2,8,3)$ & $(21395,3,971)$ \\\hline
   $(1,1,2)$ & $(3,8,2)$ & $(799,7,16060)$ \\\hline
   $(1,2,1)$ & $(3,11,2)$ & $(1719,5,71953)$ \\\hline
   $(1,2,1)$ & $(9,2,3)$ & $(13,75090,2797)$ \\\hline
   $(2,1,1)$ & $(4,5,2)$ & $(6071,959,59397119)$ \\\hline
   $(1,1,2)$ & $(2,5,4)$ & $(49800547,953,6357)$ \\\hline
   $(1,2,1)$ & $(3,7,2)$ & $(1346695,177,1566280459)$ 
\\\hline
\end{tabular}   
\end{center}
\end{table}
Conjecture \ref{conj:abcfixed} is a far-reaching generalization of Theorem \ref{th:genfermat} that allows $p,q,r$ to also be variables. There is no triple $(a,b,c)$ for which this conjecture has been proved. Conjecture \ref{conj:primesfixed} is even more general.

Eq.~\eqref{eq:fermatpqr} corresponds to the case $a=b=c=1$ of Conjecture \ref{conj:abcfixed}. The ``next'' case one may consider is $| abc| =2$. By rearranging \eqref{eq:genfermat} appropriately, we may assume that $a,b,c,x,y,z$ are all positive integers, and consider equation
\begin{equation}
\label{eq:open}
 ax^p+by^q=cz^r, \quad \text{gcd}(x,y,z)=1, \quad 1/p+1/q+1/r<1, \quad abc=2,
\end{equation}
for which we are looking for solutions in positive integers $a,b,c,p,q,r,x,y,z$.  Our computer search\footnote{This search used the ALICE High Performance Computing facility at the University of Leicester.} in the range $cz^r\leq 2^{80}$ returned only the solutions listed in 
Table \ref{tab:pqrprimesols}.  

\begin{question}
\label{q:abc2}
 Find all solutions to \eqref{eq:open} in positive integers $a$, $b$, $c$, $p$, $q$, $r$, $x$, $y$, $z$. In particular, does it have any solutions except of the ones listed in 
Table \ref{tab:pqrprimesols}?
\end{question}

While we are not sure whether the solution list to \eqref{eq:open} in 
Table \ref{tab:pqrprimesols} is complete, we conjecture that \eqref{eq:open} has no solutions with $\min(p,q,r)\geq 3$ other than $1^p+1^q=2 \cdot 1^r$, $2\cdot 5^4+3^4=11^3$ and $43^4 + 237^3=2\cdot 203^3$.

Another goal, suggested in~\cite[Section 3.4.2]{MR4620765}, is to describe \emph{all} (not necessary
primitive) integer solutions to \eqref{eq:genfermat}, or at least to \eqref{eq:fermatpqr}. This
might be easier because for some exponents $(p,q,r)$ there exist formulas describing all integer
solutions to these equations that give no hint on the structure of primitive solutions. As mentioned
above, such formulas have been derived in~\cite{MR4620765} for Eq.~\eqref{eq:genfermat} under the
condition that at least one of the integers $p,q,r$ is coprime with the other two. This condition
fails, if, for example, $(p,q,r)=(PQ,QR,RP)$ for some integers $P,Q,R$ greater than $1$. In this
case, \eqref{eq:genfermat} reduces to $ax^{PQ}+by^{QR}=cz^{RP}$, or $a+b(y^R/x^P)^Q=c(z^R/x^Q)^P$. This equation can be written as
\begin{equation}
\label{eq:cupmbvqma}
 cU^P-bV^Q=a,
\end{equation}
where $U=z^R/x^Q$ and $V=y^R/x^P$ are rational variables. If $P+Q\geq 7$, this is an equation of
genus $g\geq 2$, hence it has finitely many rational solutions by Falting's
Theorem~\cite{MR718935}, and, if we could list them all, we would be able to easily describe all
possible $x,y,z$. In the case $a=b=c=1$, \eqref{eq:cupmbvqma} reduces to $U^P-V^Q=1$.
Mih\u{a}ilescu~\cite{MR2076124}, confirming famous conjecture of Catalan, proved that the only
positive integer solutions to this equation with $P,Q\geq 2$ is $(U,V,P,Q)=(3,2,2,3)$, but solving this
equation in rational $U,V$ seems to be more difficult, see~\cite[Theorem 12.4]{MR891406}
and~\cite[Theorem 5]{mihailescu2007cylcotomic} for a partial progress.

All conjectures and open questions mentioned in this section look very difficult.\footnote{In 2022, Mochizuki et al.~\cite{mochizuki2022explicit} claimed to prove the effective abc conjecture, stating that \eqref{eq:abc} has only finitely many solutions in coprime positive integers $(A,B,C)$, and, moreover, there is an explicitly computable upper bound for the size of any such solution. If true, this would imply the truth of Conjectures \ref{conj:abcfixed} and \ref{conj:primesfixed}, the full resolution of \eqref{eq:fermatpqr} for all but finitely many explicit signatures~\cite{zhongpeng2025the}, reduce Open Question \ref{q:abc2} to finite computation, and solve many other cases of \eqref{eq:genfermat}. However, Scholze and Stix~\cite{scholze2018abc} have found a serious flaw in the argument.} Readers
preferring easier-looking open problems are invited to investigate equations listed in 
Table \ref{tab:H60Fermat}, see Section \ref{sec:system}.

\section{Equations with special signatures}
\label{sec:special}

In this section we will discuss Eq.~\eqref{eq:genfermat} with special signatures $(p,q,r)$ presented in 
Table \ref{tb:triples}. We start with the most studied case $(p,q,r)=(4,4,4)$.

\subsection{Equations of signature $(4,4,4)$}

The case $p=q=r=4$ in \eqref{eq:genfermat} results in the equation
\begin{equation}
\label{eq:ax4pby4mcz4}
 ax^4+by^4=cz^4
\end{equation}
for positive integers $a,b,c$. Sections 6.2.2 and 6.3.4 of~\cite{mainbook} and
Sections 6.4 and 6.5 of~\cite{wilcox2024systematic} study
Eqs.~\eqref{eq:ax4pby4mcz4} ordered by $a+b+c$, decide the existence of a solution $(x,y,z)\neq (0,0,0)$
for all equations with $a+b+c\leq 62$, and solves Problem \ref{prob:existence} for the equations with
$a+b+c<39$. The first Eq.~\eqref{eq:ax4pby4mcz4} for which Problem \ref{prob:existence} is left open
in~\cite{mainbook} is
\begin{equation*}\label{eq:7x4p25y4m7z4}
 7x^4+25y^4=7z^4.
\end{equation*}
However, there are Eqs.~\eqref{eq:ax4pby4mcz4} with $a+b+c<39$ that have easy-to-find non-trivial solutions, which solves Problem  \ref{prob:existence} but leaves open Problem \ref{prob:main}. For example, this happens if $a+b=c$, in which case $(x,y,z)=(\pm 1, \pm 1, \pm 1)$ are the solutions.    
The smallest non-trivial example in this category is the equation
\begin{equation*}
x^4+2y^4=3z^4,
\end{equation*}
whose only primitive solutions are $(x,y,z)=(\pm 1, \pm 1, \pm 1)$, as proved by Jeremy Rouse on Mathoverflow website.\footnote{\url{https://mathoverflow.net/questions/480114}} Other non-trivial examples with $a+b=c$ are the equations
\begin{equation*}
(a) \,\, x^4+3y^4=4z^4, \quad \text{and} \quad (b) \,\, x^4+8y^4=9z^4.
\end{equation*}
Lucas (see~\cite[pp. 630]{dickson1920history}) proved that these equations have no primitive integer
solutions other than $(\pm 1, \pm 1, \pm 1)$.

In the case $a=b=1$, Eq.~\eqref{eq:ax4pby4mcz4} reduces to
\begin{equation}
\label{eq:x4py4ecz4}
 x^4+y^4=cz^4.
\end{equation}
Problem \ref{prob:existence} for this equation reduces to the question whether a positive integer
$c$ can be represented as a \emph{sum} of two rational fourth powers. Cohen~\cite[Section
6.6]{MR2312337} solved this problem for all integers $c$ in the range $1\leq c \leq 10,000$, except of
$c=4481$, $7537$ and $8882$. These values of $c$ have been addressed
in~\cite{MR709852,MR4458887}, hence Problem \ref{prob:existence} for \eqref{eq:x4py4ecz4} is now
solved for all $1\leq c \leq 10,000$. 
In particular, $c=5906=(25/17)^4+(149/17)^4$ is the smallest integer that is the sum of two fourth powers of rational numbers but not the sum of two fourth powers of integers. 

Grigorov and Rizov~\cite{Grigorov1998Heights} proved that if $c>2$ is a fourth-power free
integer and the rank of elliptic curve $v^2=u^3-cu$ is $1$, then Eq.~\eqref{eq:x4py4ecz4} has
no non-zero integer solutions.

For values of $c$ when \eqref{eq:x4py4ecz4} is solvable, such as $c=17$, the question
to find \emph{all} primitive solutions (that is, Problem \ref{prob:main}) has not been addressed
in~\cite[Section 6.6]{MR2312337}. The mentioned case $c=17$ has been solved in 2001 by Flynn and
Wetherell~\cite{flynn2001covering}. 

\begin{theorem}
\label{th:flynn2001}
 The primitive positive integer solutions to Eq.~\eqref{eq:x4py4ecz4} with $c=17$, that is,
\begin{equation*}
 x^4+y^4=17z^4
\end{equation*}
 are $(x,y,z)=(2,1,1)$ and $(1,2,1)$.
\end{theorem}

A combination of Theorem \ref{th:flynn2001} with elementary methods discussed in the introduction solves
Eq.~\eqref{eq:x4py4ecz4} for all $1\leq c \leq 81$, and in fact for all $1\leq c\leq 100$ except of
$c=82$ and $c=97$. In 2023, using the methods developed to prove Theorem \ref{th:flynn2001} with
some extra tweaks, P{\u{a}}durariu~\cite{puadurariu2023rational} was able to resolve the case
$c=97$.

\begin{theorem}
\label{th2:flynn2001}
 The primitive positive integer solutions to Eq.~\eqref{eq:x4py4ecz4} with $c=97$, that is,
\begin{equation*}
 x^4+y^4=97z^4
\end{equation*}
 are $(x,y,z)=(2,3,1)$ and $(3,2,1)$.
\end{theorem}

To the best of our knowledge, the case $c=82$ remains open. 

Eq.~\eqref{eq:x4py4ecz4} with $z\neq 0$ can be rewritten as $(x/z)^4+(y/z)^4=c$, and reduces to finding a representation of $c$ as a sum of two rational fourth powers. A similar problem of representing an integer $b$ as a \emph{difference} of two fourth powers reduces to equation
\begin{equation}
\label{eq:x4pby4ez4}
 x^4+b y^4=z^4
\end{equation} 
with $y\neq 0$. The problem of \emph{existence} of integer solutions to \eqref{eq:x4pby4ez4} has
been solved in~\cite[Section 6.3.4]{mainbook} for $1\leq b \leq 218$, and is left open for $b=219$.
However, if a non-trivial solution to \eqref{eq:x4pby4ez4} exists, the problem of finding all
solutions has not been studied in~\cite{mainbook}. In 2020, S\"oderlund~\cite{soderlund2020note}
resolved this problem for $b=34$.

\begin{theorem}
\label{th:soderlund}
 The only primitive non-zero integer solutions to the equation
\begin{equation*}
 x^4+34 y^4=z^4
\end{equation*}
 are $(x, y, z) = (\pm 3, \pm 2, \pm 5)$.
\end{theorem}

Taclay and Bacani~\cite{taclay2023} consider Eq.~\eqref{eq:x4pby4ez4} with $b=2p$ where
$p$ is a prime and they prove the following result.

\begin{theorem}
\label{th:x4p2py4mz4}
 If $p$ is a prime satisfying any of the following conditions:
\begin{itemize}
  \item[(i)] {$p \not \equiv 1$} (mod $16$),
  \item[(ii)] {$p\equiv 3,4$} (mod $5$),
  \item[(iii)] {$p \equiv 7,8,11$} (mod $13$),
  \item[(iv)] {$p \equiv 4,5,6,9,13,22,28$} (mod $29$). 
\end{itemize}
 Then the equation
\begin{equation}
\label{eq:x4p2py4mz4}
  x^4+2py^4=z^4
\end{equation}
 has no non-trivial primitive integer solutions.
\end{theorem}

The smallest prime $p$ for which Eq.~\eqref{eq:x4p2py4mz4} is not covered by Theorems \ref{th:soderlund} and \ref{th:x4p2py4mz4} is $p=97$. 

Eq.~\eqref{eq:ax4pby4mcz4} with coefficients $a,b,c$ being perfect cubes was studied
in~\cite{MR249355,Mordell_1970}, and we have the following result. 
\begin{theorem}
\label{th:MR249355}
 Let $s,t,u$ be pairwise coprime positive integers, with $s \equiv t \equiv u \equiv -1$ (mod 8). Let $n$ be an integer divisible by $8$ and let $v,w$ be integers satisfying $v \equiv w \equiv -1$ (mod 8), and ratios
\begin{equation*}
 \frac{tw-uv}{s}, \quad \frac{un-sw}{t}, \quad \frac{sv-tn}{u}
\end{equation*}
 are integers whose positive odd factors are all $1$ modulo $8$. Then, there are no non-trivial integer solutions of
\begin{equation*}
 \left(\frac{tw-uv}{s} \right)^3 x^4 + \left(\frac{un-sw}{t} \right)^3 y^4 = \left(\frac{tn-sv}{u} \right)^3z^4.
\end{equation*}
\end{theorem}
Examples of values $(s,t,u,n,v,w)$ satisfying these conditions are 
$(7,15,23,\\ 8280,4991,13335)$, $(7,15,23,8280,16583,15855)$ and $(7,15,23,11040,3703, \\14175)$ which correspond to equations
\begin{equation*}
\begin{aligned}
& 12176^3x^4 + 6473^3 y^4=3881^3z^4, \quad 20512^3x^4 + 353^3 y^4=5297^3z^4 , \\ & 18208^3x^4 + 10313^3 y^4=6073^3z^4 ,
\end{aligned}
\end{equation*} 
respectively.

\subsection{Equations of signature $(2,4,r)$}

The case $p=2$, $q=4$ and prime $r \geq 5$ in \eqref{eq:genfermat} results in the equation
\begin{equation}
\label{eq:sig24r}
 ax^2+by^4=cz^r, \quad r \geq 5
\end{equation}
with positive integer parameters $a,b,c$. 

Let us first consider Eq.~\eqref{eq:sig24r} with $a=2$ and $b=c=1$, that is,
\begin{equation}
\label{eq:x4p2y2mzr}
 2x^2+y^4 = z^r, \quad\quad r\geq 5.
\end{equation}
In 2008, Dieulefait and Urroz~\cite{dieulefait2008solvingfermattypeequations} proved that if
$r >349$ is a prime, then Eq.~\eqref{eq:x4p2y2mzr} does not have non-trivial primitive
solutions. In 2010, Bennett, Ellenberg and Ng~\cite{bennett2010diophantine} significantly strengthened
this result and completely solved Eq.~\eqref{eq:x4p2y2mzr} for all integers $r \geq 5$. 

\begin{theorem}
\label{th:bennett2010rint}
  Eq.~\eqref{eq:x4p2y2mzr} has no non-trivial primitive solutions for integer $r\geq 6$, while for $r=5$ its only non-trivial primitive integer solutions are $(x,y,z)=(\pm 11, \pm 1,3)$.  
\end{theorem}

Let us now consider Eq.~\eqref{eq:sig24r} with $a=3$ and $b=c=1$. Dieulefait and
Urroz~\cite{dieulefait2008solvingfermattypeequations} proved that if $r >131$ is a prime, then the
equation
\begin{equation*}
3x^2+y^4= z^r,
\end{equation*}
does not have any non-trivial primitive solutions. Pacetti and Villagra Torcomian~\cite{MR4473105}
remarked that this result can be improved to $r > 17$ by Proposition 5.4
of~\cite{koutsianas2019generalizedfermatequationa23b6cn}.

In 2022, Pacetti and Villagra Torcomian~\cite{MR4473105} studied Eq.~\eqref{eq:sig24r} with
$a \in \{5,6,7\}$, $b=c=1$ and prime $p$, that is,
\begin{equation}
\label{eq:ax4p5y2pzp}
 ax^2 + y^4 = z^p,
\end{equation}
and they proved the following results.

\begin{theorem}
\label{th:ax4p5y2pzp}
   Let $(a,p_0)\in\{(5,499), (6,563), (7,349)\}$. 
 Assume that $p > p_0$ is prime. Then Eq.~\eqref{eq:ax4p5y2pzp} has no non-trivial primitive solutions.
 
\end{theorem}

D\c abrowski~\cite{MR2737959} studied equations of the form $y^4- z^2 = sx^p$, which, after substitution
$x\to -x$, can be rewritten with positive coefficients as
\begin{equation}
\label{eq:qxppy4z2}
 sx^p + y^4 = z^2,
\end{equation}
where $(s,p)$ are certain pairs of odd primes. They used a result from~\cite{ivorra2006quelques}
to deduce the following theorem.

\begin{theorem}
\label{th:qxppy4z2}
 Let $s > 3$ be a prime and let $s \equiv 3$ mod 8 and $s \neq 2t^2 + 1$,
 or $s \equiv 5$ mod 8 and $s \neq t^2 + 4$. In addition, let $p$ be a prime satisfying
 $p > (8 \sqrt{s+1}+1)^{16(s-1)}$. Then Eq.~\eqref{eq:qxppy4z2}  has no non-trivial primitive integer solutions.
\end{theorem}

Examples of $s$ satisfying the conditions of Theorem \ref{th:qxppy4z2} are $s=11$, $37$, $43$ and so on.

Equation
\begin{equation}
\label{eq:qxppy4z4}
 sx^p + y^4 = z^4
\end{equation}
is a special case of \eqref{eq:qxppy4z2} with $z$ in \eqref{eq:qxppy4z2} being a perfect
square. This equation has been studied in~\cite{Dabrowski_2007,MR2737959}, where the following results
have been obtained. 
\begin{theorem}
Let $\alpha\geq 0$ be an integer and $p\geq 5$ be a prime. Then equation
\begin{equation*}
 2^\alpha x^p + y^4 = z^4
\end{equation*}
 has no non-trivial primitive integer solutions.
\end{theorem} 

\begin{theorem}
\label{th:2asbxppy4mz4}
 Let $s$ be an odd prime and $s \neq 2^t \pm 1$. Let $p$ be a prime satisfying $p > (\sqrt{
  8s + 8} + 1)^{2s-2}$. Let $\alpha\geq 0$ and $\beta>0$ be integers. Then equation
\begin{equation}
\label{eq:2asbxppy4mz4}
  2^\alpha s^\beta x^p + y^4 = z^4
\end{equation}
 has no non-trivial primitive integer solutions.
\end{theorem} 

Eq.~\eqref{eq:2asbxppy4mz4} has been further studied by Bennett~\cite{MR4205757}, who proved a
version of Theorem \ref{th:2asbxppy4mz4} with a smaller set of excluded primes.

\begin{theorem}
\label{th:bennett21}
 Let $s$ be a prime not of the form $s=2^{2^k}+1$ for integer $k\geq 1$, and let $\alpha,\beta$ be non-negative integers. Let $p$ be a prime satisfying $p > (\sqrt{8s + 8} + 1)^{2s-2}$. Then Eq.~\eqref{eq:2asbxppy4mz4} has no non-trivial primitive integer solutions.
\end{theorem} 

Primes of the form $s=2^{2^k}+1$ are called Fermat primes. It is widely believed that the only such
primes are $s=3,5,17,257$ and $65537$. The case $s=3$ corresponds to $k=0$, while Theorem
\ref{th:bennett21} fails to apply only if $k\geq 1$. Bennett~\cite{MR4205757} also solved Eq.~\eqref {eq:2asbxppy4mz4} for the exceptional cases $s=5$ and $17$, at least for
$(\alpha,\beta)=(0,1)$.

\begin{theorem}
For every prime $p>5$, equations
\begin{equation*}
 5 x^p + y^4 = z^4 \quad\quad \text{and} \quad\quad 17 x^p + y^4 = z^4 
\end{equation*}
 have no non-trivial primitive integer solutions.
\end{theorem} 

Theorem \ref{th:bennett21} is applicable only if $p$ is large in terms of $s$. For small
values of $s$, Bennett~\cite{MR4205757} was able to remove this condition and extend Theorem
\ref{th:bennett21} to all $p>5$.

\begin{theorem}
If $s$ is a prime with $2\leq s < 50$, $s\neq 5$, $s\neq 17$, $\alpha$ and $\beta$ are non-negative integers, and $p>5$ is a prime, then Eq.~\eqref{eq:2asbxppy4mz4} has no non-trivial primitive integer solutions.
\end{theorem} 

Pacetti and Villagra Torcomian~\cite{MR4609012} studied equation
\begin{equation}
\label{eq:xppby2mz4}
 x^p+by^2=z^4
\end{equation}
for certain values of $b$, and they obtained the following results. 

\begin{theorem}
Let $p > 19$ be a prime number such that $p \neq 97$ and $p \equiv 1, 3$ modulo $8$. Then, Eq.~\eqref{eq:xppby2mz4} with $b=6$
       has the only non-trivial primitive integer solutions $(x,y,z)=(1, \pm 20, \pm 7)$.
\end{theorem}

\begin{theorem}
Let $p > 19$ be a prime number such that $p \equiv 1, 3$ modulo $8$. If
\begin{itemize}
  \item[(i)] {$b=10$} and $p \neq 139$, or
  \item[(ii)] {$b=11$} and $p \neq 73$, or
  \item[(iii)] {$b=19$} and $p \neq 43,113$, or
  \item[(iv)] {$b=129$} and $p \neq 43$,
\end{itemize}
 then Eq.~\eqref{eq:xppby2mz4} has no non-trivial primitive solutions.
\end{theorem}

\begin{theorem}
Let $p > 900$ be a prime number. Then, there are no non-trivial primitive solutions of Eq.~\eqref{eq:xppby2mz4} with $b=129$.
\end{theorem}

Cao~\cite{MR1341665} studied Eq.~\eqref{eq:xppby2mz4} with $b=p$, that is,
\begin{equation}
\label{eq:xpppy2mz4}
 x^p+py^2=z^4
\end{equation}
and proved that if $p \equiv 1$ modulo $4$ is a prime and\footnote{Here, $B_n$ denotes the $n$th Bernoulli number.} $p\nmid B_{(p-1)/2}$, then Eq.~\eqref{eq:xpppy2mz4} has no integer solutions with $\gcd(y,z)=1$, $p\vert y$ and $2\vert z$.

Langmann~\cite{MR1604052} considers the equation
\begin{equation}
\label{eq:ax2mpy2mzn}
 ax^{2m}+y^2=z^n
\end{equation}
and they obtained the following result.

\begin{theorem}
\label{th:ax2mpy2mzn}
 Let $m,n \geq 2$ be integers. Then, for almost all square-free $a \not \equiv -1$ modulo $4$ with
\begin{equation*}
 \frac{n}{\gcd(n,h(-a))} \geq \max\{6,14-2m\},
\end{equation*}
 where $h(-a)$ denotes the class number\footnote{A list of the class numbers for the first 10,000 square-free $a$ can be found at \url{https://oeis.org/A000924/b000924.txt}.} of $\mathbb{Q}(\sqrt{-a})$, Eq.~\eqref{eq:ax2mpy2mzn} has no non-trivial primitive integer solutions. 
\end{theorem}
The case $m=2$ of Theorem \ref{th:ax2mpy2mzn} solves many equations of the form $ax^4+y^2=z^n$ of signature $(4,2,n)$.

\subsection{Equations of signature $(2,6,r)$}

Some researchers have considered equations with signatures $(2,6,r)$ for $r\geq 3$. This family contains special signatures $(2,6,4)$, $(2,6,5)$ and $(2,6,6)$.

Chen~\cite{Chen_2012} studied Eq.~\eqref{eq:genfermat} with $p$ prime, $q=6$,
$b=r=2$ and $a=c=1$, that is,
\begin{equation}
\label{chen_a1c1}
 x^p +2y^6= z^2.
\end{equation}

\begin{theorem}
Let $p$ be a prime such that $p \equiv 1, 7$ (mod $24$) and $p \neq 7$. Then Eq.~\eqref{chen_a1c1} does not have any non-trivial primitive integer solutions except those with $x = \pm1$.
\end{theorem}

Eq.~\eqref{eq:genfermat} with $r$ prime, $q=6$, $p=2$ and $a=c=1$, that is,
\begin{equation}
\label{eq:x2pby6mzr}
 x^2+by^6=z^r,
\end{equation}
has been studied in several papers. The case $b=2$ of \eqref{eq:x2pby6mzr}, that is,
\begin{equation}
\label{eq:x2p2y6mzr}
 x^2 + 2y^6 = z^r,
\end{equation}
has been studied by Pacetti and Villagra Torcomian~\cite{MR4473105}, where it was resolved for all
$r>257$.

\begin{theorem}
Let $r > 257$ be a prime number. Then Eq.~\eqref{eq:x2p2y6mzr} has no non-trivial primitive integer solutions.
\end{theorem}

Koutsianas~\cite{koutsianas2019generalizedfermatequationa23b6cn} solved the case $b=3$ of
\eqref{eq:x2pby6mzr}, that is, equation
\begin{equation}
\label{eq:x2p3y6mzr}
 x^2 +3y^6 = z^r.
\end{equation}

\begin{theorem}
Let $r \geq 3$ be an integer. If $r\neq 4$, then Eq.~\eqref{eq:x2p3y6mzr} has no non-trivial primitive integer solutions, while for $r=4$ its only non-trivial primitive integer solutions are $(x, y, z) = (\pm 47, \pm 2, \pm 7)$. 
\end{theorem}

Eq.~\eqref{eq:x2pby6mzr} with $b=6$, that is,
\begin{equation}
\label{eq:x2p6y6mzr}
 x^2 + 6y^6 = z^r,
\end{equation}
was solved by Pacetti and Villagra Torcomian~\cite{MR4473105} for $r>563$.
\begin{theorem}
Let $r > 563$ be a prime number. Then, Eq.~\eqref{eq:x2p6y6mzr} has no non-trivial primitive integer solutions.
\end{theorem}

Eq.~\eqref{eq:x2pby6mzr} with some other values of $b$ has been studied in~\cite{MR4583916},
where the following theorem has been proven. 

\begin{theorem}
\label{th:x2pby6mzr}
 Let $(b,r_0)$ be as in 
Table \ref{tb:x2pby6mzr}. Assume that $r\geq r_0$ is a prime satisfying the congruence conditions in 
Table \ref{tb:x2pby6mzr}. Then Eq.~\eqref{eq:x2pby6mzr} has no non-trivial primitive integer solutions.
\end{theorem}

\begin{table}
\begin{center}
  \caption{\label{tb:x2pby6mzr}Details for Theorem \ref{th:x2pby6mzr}. }
\begin{tabular}{|c|c|c|}
\hline

   $b$ & $r_0$ & Congruence conditions \\   \hline\hline
   5 & 1033 & \\\hline
   7 & 337 & $r \equiv 5, 7$ (mod 12) \\\hline
   11 & 557 & $r \equiv 3$ (mod 4) \\\hline
   13 & 3491 & \\\hline
   15 & 743 & $r \equiv 5, 7, 15, 17, $ or 19 (mod 24) \\\hline
   19 & 1031 & 
\\\hline
\end{tabular}
\end{center}
\end{table}

In the same paper~\cite{MR4583916}, the authors remarked that they also studied \eqref{eq:x2pby6mzr}
with $b=10$ and $b=17$, but the resulting computation is unfeasible.

Zelator~\cite{MR1188732} proves that for any positive odd integer $b$ with a prime factor
$l \equiv \pm 3$ modulo $8$ and positive integers $m$ and $n$, the equation
\begin{equation*}
x^2 + b^2y^{2m} = z^{4n}
\end{equation*}
has no positive integer solutions with $\gcd(x,by)=1$ and $y$ odd. With $m=3$ and $n=1$, this gives some partial progress towards the resolution of the equation $x^2 + b^2y^{6} = z^{4}$ of signature $(2,6,4)$. 

\subsection{Equations of signature $(4,4,3)$}

Let us now discuss Eq.~\eqref{eq:genfermat} with $(p,q,r)=(4,4,3)$ and $a=c=1$, that is, equation
\begin{equation}
\label{eq:x2pby4pz3}
 x^4+b y^4=z^3.
\end{equation}
The case $b=2$ of \eqref{eq:x2pby4pz3} has been resolved in 2017 by
S\"oderlund~\cite{soderlund2017primitive}. 

\begin{theorem}
\label{th:x4p2y4pz3}
 The only non-trivial primitive integer solutions to the equation
\begin{equation*}
 x^4+2y^4=z^3
\end{equation*}
 are $(x,y,z)=(\pm 3, \pm 5, 11)$.
\end{theorem}

The case $b=3$ was solved in 1994 by Terai~\cite{MR1288426}.
\begin{theorem}
The equation
\begin{equation*}
 x^4+3y^4=z^3
\end{equation*}
 has no non-trivial primitive integer solutions.
\end{theorem}

The next case $b=4$ is a special case of the following theorem~\cite{soderlund2020diophantine}.

\begin{theorem}
For any prime number $s$, equation
\begin{equation*}
 x^4+s^2 y^4=z^3
\end{equation*}
 has no non-trivial primitive integer solutions.
\end{theorem}

The next case is $b=5$, and we do not know whether this case has been solved. In 2019,
S\"oderlund~\cite{soderlund2019some} reproved the case $b=3$ and resolved the cases of
$b=19,43$ and $67$.

\begin{theorem}
For $b \in \{3,19,43,67\}$, equation
\begin{equation*}
 x^4+b y^4=z^3
\end{equation*}
 has no non-trivial primitive integer solutions.
\end{theorem}

\subsection{A systematic approach}
\label{sec:system}

\begin{table}
\begin{center}
  \caption{\label{tab:H60Fermat}Uninvestigated equations \eqref{eq:genfermat} with $H=| a| 2^p+| b| 2^q+| c| 2^r \leq 60$.}
\begin{tabular}{|c|c|c|c|c|c|}
\hline
$H$ & Equation & $H$ & Equation & $H$ & Equation \\
\hline \hline
   40 & $x^4+2y^3+z^3=0$ & 56 & $x^4+4y^3+z^3=0$  &56 & $x^5+y^4-2z^2=0$     \\ 
   \hline
   48 & $x^4+3y^3+z^3=0$ & 56 & $2x^4-y^4+z^3=0$  & 60 & $x^5+y^4+3z^2=0$    \\ 
   \hline
   56 & $x^4+3y^3+2z^3=0$  &56 & $x^5+2y^3+z^3=0$   &60 & $x^5+y^4-3z^2=0$ 
\\\hline
\end{tabular}
\end{center}
\end{table}

Monograph~\cite{mainbook} suggests studying Diophantine equations systematically. It defines the size
of Eq.~\eqref{eq:genfermat} as $H=| a| 2^p+| b| 2^q+| c| 2^r$, and lists all non-trivial equations of size $H\leq 60$.
Equations of small size with $| abc| =1$ are covered in 
Table \ref{tb:pqrsigs}. Equations $x^4+y^4=2z^3$ and $y^3+z^3=2x^4$ of size $H=48$ are special cases of Theorems
\ref{th2:bennett2004} and \ref{th:zhang2014}, respectively.  Equation $2x^2+y^4=z^5$ of size $H=56$ is the
$r=5$ case of Theorem \ref{th:bennett2010rint}. Equation $x^4+2y^4=z^3$ of the same size has been
solved in Theorem \ref{th:x4p2y4pz3}. Equations $x^4+2y^3+2z^3=0$, $x^4-y^4=2z^3$ and $x^4-y^4=3z^3$ of sizes
$H=48$, $48$ and $56$, respectively, have been solved in~\cite[Appendix
A]{wilcox2024systematic}. All other non-trivial equations of size $H\leq 60$ are listed in 
Table \ref{tab:H60Fermat}. All these equations except $x^5+2y^3+z^3=0$ have special signatures, and we were not able to find a reference where any of them has been solved. We invite the reader to investigate these equations.

\section{Equations of signature $(p,p,p)$}

\subsection{Reduction to computing rational points on hyperelliptic curves}
\label{sec3.1}

As discussed in the introduction, if the exponents $(p,q,r)$ in \eqref{eq:genfermat} are not special, then we may assume that all $p,q,r$ are prime. In this section we discuss the case $p=q=r$, that is, equation
\begin{equation}
\label{eq:axppbypmczp}
 ax^p+by^p=cz^p,
\end{equation}
where $a,b,c$ are positive integers, and $p\geq 5$ is a prime. It is well-known~\cite[Section
6.2.3]{mainbook} and easy to check that if $(x,y,z)$ with $xyz\neq 0$ is an integer solution to
\eqref{eq:axppbypmczp}, then
\begin{equation*}
\begin{aligned}
 (X,Y)=\,&(xy/z^2,a(x/z)^p-c/2), \quad (xz/y^2,-c(z/y)^p+b/2)\nonumber\\
 &\text{and} \quad (yz/x^2,b(y/x)^p+a/2)
\end{aligned}
\end{equation*}
are rational points on the curves
\begin{equation}
\label{eq:homdia3dhyp3}
 Y^2 = -abX^p+c^2/4, \quad Y^2 = acX^p+b^2/4 \quad \text{and} \quad Y^2 = bcX^p+a^2/4,
\end{equation}
respectively. In all cases, $X\neq 0$.
These are curves of genus $g=(p-1)/2 \geq 2$, and, by Faltings' theorem~\cite{MR718935}, each of them has a
finite number of rational points. If these points can be computed for at least one of these curves, it
is straightforward to convert them into solutions to \eqref{eq:axppbypmczp}. However, computing
rational points on a high-genus curve is a difficult problem in general, and the algorithms are known
only in some special cases. One example is given by the following theorem, which is a corollary of the
main result of~\cite{em/1069786350}, and formulated explicitly in~\cite[Theorem 3.76]{mainbook}.

\begin{theorem}
\label{th:jacrank0}
 There is an algorithm for computing all rational solutions to any equation defining a curve of genus $g\geq 2$ whose Jacobian has rank $r=0$.
\end{theorem} 

The definition of the rank of Jacobian $r$ in Theorem \ref{th:jacrank0} is too technical to be presented here, but we mention that for the equation $y^2=P(x)$, the rank $r$ can be computed (most often estimated) using the following Magma commands
\begin{align*}
	&{\tt > P<x> := PolynomialRing(Rationals());}\\
	&{\tt > C := HyperellipticCurve(P(x));}\\
	&{\tt > J := Jacobian(C);}\\
	&{\tt > RankBounds(J);}
\end{align*}
that can be run in the Magma calculator\index{Magma calculator} \url{http://magma.maths.usyd.edu.au/calc/}. In general, the output looks like ``$a\,\,b$'', where $a$ and $b$ are integers that are the lower and upper bounds for the rank $r$ of the Jacobian. If for at least one out of the three curves \eqref{eq:homdia3dhyp3} the output is ``$0\,\,0$'', then its rational points can be computed by Theorem \ref{th:jacrank0}, which in turn leads to the resolution of the corresponding Eq.~\eqref{eq:axppbypmczp}.

\subsection{Explicit small values of $p$}

We next discuss Eq.~\eqref{eq:axppbypmczp} for explicit small $p\geq 5$. In this section, we will not assume that $p$ is prime.

With $p=5$, Eq.~\eqref{eq:axppbypmczp} becomes
\begin{equation}
\label{eq:ax5pby5mcz5}
 ax^5+by^5=cz^5,
\end{equation}
where $a,b,c$ are positive integers. For Eq.~\eqref{eq:ax5pby5mcz5}, the corresponding curves
\eqref{eq:homdia3dhyp3} have genus $g=2$, and if at least one of these curves has the rank of
the Jacobian $r\leq 1$, then its rational points can be computed by Magma's {\tt Chabauty} command,
which leads to the complete resolution of \eqref{eq:ax5pby5mcz5}. As noted in
Sections 6.2.2 and 6.3.4 of~\cite{mainbook}, this method resolves Problem \ref{prob:existence} for all Eqs.~\eqref{eq:ax5pby5mcz5} with $a+b+c<19$. If
Eqs.~\eqref{eq:ax5pby5mcz5} are ordered by $a+b+c$, then the first Eq.~\eqref{eq:ax5pby5mcz5} for which Problem \ref{prob:existence} is left open in~\cite{mainbook} is
\begin{equation}
\label{eq:4x5p4y5m11z5}
 4x^5+4y^5=11z^5.
\end{equation}

If an Eq.~\eqref{eq:ax5pby5mcz5} has some obvious non-trivial integer solution, this immediately solves Problem \ref{prob:existence} with the ``Yes'' answer, but leaves open the problem of listing all primitive integer solutions (Problem \ref{prob:main}). In particular, equation
\begin{equation*}
2x^5+3y^5=5z^5
\end{equation*}
has solution $(x,y,z)=(1,1,1)$, but we do not know whether it has any other non-trivial primitive integer solutions.

Problem \ref{prob:existence} for Eq.~\eqref{eq:ax5pby5mcz5} with $a=b=1$, that is,
\begin{equation}
\label{eq:x5py5ecz5}
 x^5+y^5=cz^5,
\end{equation}
reduces to investigating which positive integers $c$ can be represented as a sum of two
rational fifth powers. As remarked in~\cite[Section 6.3.4]{mainbook}, the discussed argument solves
this question for all positive integers $1\leq c\leq 100$ except of $c = 68, 74, 87$ and $88$.

Dirichlet (see~\cite[pp. 735]{dickson1920history}) proved the non-existence of non-zero integer
solutions to \eqref{eq:x5py5ecz5} for certain values of $c$ with no prime factors congruent
to $1$ modulo $5$. Lebesgue conjectured that \eqref{eq:x5py5ecz5} is not solvable
for any such $c$, except of solutions with $x=y=z$ for $c=2$. In 2004,
Halberstadt and Kraus~\cite{kraus2004} confirmed this conjecture.

\begin{theorem}
\label{th:halberstadt2004}
 Let $c$ be an integer with no prime factors congruent to $1$ modulo $5$. If $c \neq 2$, Eq.~\eqref{eq:x5py5ecz5} has no non-trivial integer solutions. If $c =2$, the only non-trivial primitive solutions are $(x,y,z)=\pm(1,1,1)$.
\end{theorem}

In particular, out of the exceptional values $c = 68, 74, 87$ and $88$ listed above, Theorem \ref{th:halberstadt2004} covers the first three, hence the only $1\leq c \leq 100$ for which Problem \ref{prob:existence} for Eq.~\eqref{eq:x5py5ecz5} remains open is $c=88$. In this case, \eqref{eq:x5py5ecz5} reduces to \eqref{eq:4x5p4y5m11z5}.

With $p=6$, Eq.~\eqref{eq:axppbypmczp} becomes
\begin{equation}
\label{eq:ax6pby6mcz6}
 ax^6+by^6=cz^6.
\end{equation}
In 2023, Newton and Rouse~\cite{MR4552507} investigated the solvability of Eq.~\eqref{eq:ax6pby6mcz6} with
$a=b=1$. Specifically, they proved that $164634913$ is the smallest positive integer that is the
sum of two sixth powers of rational numbers, and not the sum of two sixth powers of integers. To prove
this result, the authors determine all positive integers $c\leq 164634913$ for which the equation
\begin{equation*}
x^6+y^6=c z^6
\end{equation*}
has a solution with $z\neq 0$. 

Dirichlet (see~\cite[pp. 736]{dickson1920history}) proved that Eq.~\eqref{eq:axppbypmczp} with
$p=14$, $a=c=1$ and $b=2^{\alpha} 7^{\beta}$ for $\alpha \geq 0$ and $\beta \geq 1$, that is,
\begin{equation*}
x^{14}+2^{\alpha} 7^{\beta} y^{14}=z^{14},
\end{equation*}
has no non-trivial primitive integer solutions.

\subsection{The case $a=b=1$}

We next discuss the works investigating infinite families of Eq.~\eqref{eq:axppbypmczp} in which some of $a,b,c$ are fixed while $p$ varies. We start with the case $a=b=1$. In this instance, Eq.~\eqref{eq:axppbypmczp} reduces to
\begin{equation}
\label{eq:genfermata1b1ppp}
 x^p+y^p=cz^p,
\end{equation}
where $c$ is a positive integer and $p\geq 5$ is a prime.

Eq.~\eqref{eq:genfermata1b1ppp} with $c=15$, that is,
\begin{equation}
\label{eq:xppypm15zp}
 x^p+y^p=15z^p
\end{equation}
was considered by Kraus~\cite{kraus1996equations} who proved the following result.
\begin{theorem}
\label{th:xppypm15zp}
 If $p\geq 5$ is a prime\footnote{Kraus~\cite{kraus1996equations} states the result for
$p\geq 7$, but direct verification using the method described in  \ref{sec3.1} shows that it
remains correct for $p=5$.} such that either $2p+1$ or $4p+1$ is also prime then
Eq.~\eqref{eq:xppypm15zp} has no non-trivial integer solutions.
\end{theorem}
Examples of primes satisfying the condition of Theorem \ref{th:xppypm15zp} are $p=5$, $7$,
$11$, $13$, $23$ and so on. Kraus~\cite{kraus1996equations} also proved that
for \emph{all} primes $p\geq 5$ Eq.~\eqref{eq:xppypm15zp} has no non-trivial integer
solutions with $xyz$ divisible by $p$.

We next consider Eq.~\eqref{eq:genfermata1b1ppp} with $c=p$, that is,
\begin{equation}
\label{eq:genfermata1b1cppp}
 x^p+y^p=pz^p.
\end{equation}

Maillet (see~\cite[pp. 759]{dickson1920history}) solved Eq.~\eqref{eq:genfermata1b1cppp} with
$p$ being  a regular prime. Recall that a prime $p$ is said to be a regular prime,
when it is not a factor of the numerator of any one of the first $(p-3)/2$ Bernoulli numbers.

\begin{theorem}
\label{th:regular}
 For any regular prime $p$, Eq.~\eqref{eq:genfermata1b1cppp} has no non-trivial integer solutions.
\end{theorem}

A note of Westlund~\cite{MR1517149} extended the theory of this equation for any odd prime
$p$, resulting in the following partial result.

\begin{theorem}
For any odd prime $p$, Eq.~\eqref{eq:genfermata1b1cppp} has no integer solutions such that  $z$ is coprime to $p$.
\end{theorem}

Ribet~\cite{ribet1997equation} considered Eq.~\eqref{eq:genfermata1b1ppp} with $c=2^{\alpha}$ where
$\alpha$ is a positive integer, that is,
\begin{equation}
\label{eq:ribetfermat}
 x^p+y^p =2^{\alpha} z^p.
\end{equation}
The case when $\alpha=1$ corresponds to the three integers $x^p$, $z^p$ and
$y^p$ forming an arithmetic progression. Confirming conjecture of D\'enes'~\cite{MR68560},
Darmon and Merel~\cite{Darmon1997,MR1730439} proved that this is impossible unless $xyz=0$ or
$x$, $y$ and $z$ are all equal.

\begin{theorem}
\label{th:Darmon1997}
  Eq.~\eqref{eq:genfermata1b1ppp} with $c=2$ has no primitive solutions with $| xyz| >1$ for integers $p \geq 3$.
\end{theorem}

Theorem \ref{th:Darmon1997} treats the case of $\alpha=1$ in \eqref{eq:ribetfermat}. For $2\leq \alpha <p$,
Ribet~\cite{ribet1997equation} proved the following result.

\begin{theorem}
\label{th2:Darmon1997}
 For prime $p\geq 3$ and any $2\leq \alpha <p$, Eq.~\eqref{eq:ribetfermat} has no solutions with $xyz \neq 0$.
\end{theorem}

Because the case of general $\alpha\geq 0$ can be reduced to the case $0\leq \alpha <p$, and the case $\alpha=0$ is covered by Fermat's Last Theorem, a combination of Theorems \ref{th:Darmon1997} and \ref{th2:Darmon1997} solves \eqref{eq:ribetfermat} for all $\alpha\geq 0$.

We next discuss Eq.~\eqref{eq:genfermata1b1ppp} with $c=s^{\alpha}$ for odd prime $s$, that is, equation
\begin{equation}
\label{eq:ribetfermats}
 x^p + y^p = s^{\alpha} z^p.
\end{equation}

Kraus~\cite{MR1611640} investigated Eq.~\eqref{eq:ribetfermats} and proved the following result.

\begin{theorem}
\label{th:kraus1997}
 Let $s$ be an odd prime number which is not of the form $s=2^k\pm 1$, and let $\alpha$ be a non-negative integer. If $p$ is a prime such that $p > \left(1+\sqrt{(s+1)/2}\right)^{(s+11)/6}$, then Eq.~\eqref{eq:ribetfermats} has no solutions in non-zero integers $x,y,z$. 
\end{theorem}

Theorem \ref{th:kraus1997} works only if $p$ is large in terms of $s$. For small
$s$, we have the general results covering all $p$. As remarked by
Ribet~\cite{ribet1997equation}, the following result was proved by Serre~\cite{MR885783} subject to
some conjectures that have been established in later
works~\cite{MR1611640,MR1047143,wiles1995modular}.

\begin{theorem}
\label{th:ribet1997}
 Suppose that $p$ and $s$ are prime numbers, $3 \leq s \leq 60$, $s\neq 31$, $p \geq 11$, $p\neq s$. 
 If $\alpha \geq 1$, then there are no triples of non-zero integers $(x, y, z)$ which satisfy \eqref{eq:ribetfermats}.
\end{theorem}

Cohen~\cite[Theorem 15.5.3]{MR2312338} extended Theorem \ref{th:ribet1997} to the range $3 \leq s \leq 100$ and all
$p\geq 5$. 

\begin{theorem}
\label{th:ribetfermats}
 Suppose that $3\leq s \leq 100$ is a prime, $s\neq 31$. Then Eq.~\eqref{eq:ribetfermats} with prime $p \geq 5$ and any $\alpha\geq 1$, does not have any solutions with $x,y,z$ non-zero and pairwise coprime. 
\end{theorem}

Theorem \ref{th:ribetfermats} does not work for $s=31$, because the resulting equation
\begin{equation}
\label{eq:xppypm31alphazp}
 x^p + y^p = 31^{\alpha} z^p,
\end{equation}
has the primitive solution $(x, y, z) = (2, -1, 1)$ when $\alpha =1$ and $p=5$. 
For $7 \leq p \leq 10^6$, Cohen~\cite{MR2312338} proved the following theorem. 
\begin{theorem}
\label{th:xppypm31alphazp}
 Eq.~\eqref{eq:xppypm31alphazp} does not have any non-trivial pairwise coprime solutions $(x,y,z)$ for prime $p$ satisfying
 $7 \leq p \leq 10^6$.
\end{theorem}
See~\cite[Theorem 15.6.4]{MR2312338} for the case $11\leq p \leq 10^6$ and~\cite[Corollary 15.7.5]{MR2312338}
for the case $p=7$.

We next return to Eq.~\eqref{eq:genfermata1b1ppp} with $c$ not necessarily being a prime power.
Sitaraman~\cite{sitaraman2000fermat} considered the case $c=sp$ for regular prime $p$
and integer $s$, and obtained the following theorem. 

\begin{theorem}
Let $p\geq 5$ be a regular prime, and let $s$ be an integer divisible only by primes of the form $kp-1$ where $\gcd(k, p)=1$. Then
\begin{equation*}
 x^p+ y^p= ps z^p
\end{equation*}
 has no non-trivial integer solutions.
\end{theorem}

\subsection{The case $a=1$}

In this instance, Eq.~\eqref{eq:axppbypmczp} reduces to
\begin{equation}
\label{eq:xppbyppczp}
 x^p+by^p=cz^p,
\end{equation}
where integers $b,c>1$.

Eq.~\eqref{eq:xppbyppczp} with $b=3$, $c=5$ and $p \geq 7$ a prime, that is,
\begin{equation}
\label{eq:xpp3ypp5zp}
 x^p+3y^p=5z^p
\end{equation}
was considered by Kraus~\cite{kraus1996equations} who proved that the following result.

\begin{theorem}
\label{th:xpp3ypp5zp}
 If $p \geq 5$ is a prime\footnote{Kraus~\cite{kraus1996equations} states the result for
$p \geq 7$, but direct verification shows that it remains correct for $p=5$.} such that
$2p+1$ or $4p+1$ is also prime then Eq.~\eqref{eq:xpp3ypp5zp} has no non-zero integer
solutions.
\end{theorem}
Examples of primes satisfying the condition of Theorem \ref{th:xpp3ypp5zp} are $p=5$, $7$,
$11$, $13$, $23$ and so on. Kraus~\cite{kraus1996equations} also proved that
for \emph{all} primes $p\geq 5$ Eq.~\eqref{eq:xpp3ypp5zp} has no non-trivial integer
solutions with $xyz$ divisible by $p$. Kraus~\cite{Kraus2002} later verified that for
$5 \leq p <10^7$, Eq.~\eqref{eq:xpp3ypp5zp} has at least one local obstruction and therefore has no
non-trivial integer solutions.

Kraus~\cite{MR1611640} considered equations of the form \eqref{eq:xppbyppczp} with coefficients being
powers of primes, and obtained the following results.
\begin{theorem}
\label{th1:MR1611640} 
 Let $s$ be an odd prime number not of the form $s=2^k\pm 1$. Let $\alpha$ and $\beta$ be non-negative integers. If $p$ is a prime such that $p > \left(1+\sqrt{8(s+1)}\right)^{2(s-1)}$, $b=2^{\beta}$, $c=s^{\alpha}$ and if $\alpha$ is not a multiple of $p$, then there are no integer solutions to Eq.~\eqref{eq:xppbyppczp} with $xyz \neq 0$.
\end{theorem}

\begin{theorem}
\label{th2:MR1611640}
 Let $s$ be an odd prime number which is not equal to $17$, and $p \geq 5$ is prime with $p \neq s$. If $p > \left(1+\sqrt{(s+1)/6}\right)^{(s+1)/6}$, $b=16$ and $c=s^{\alpha}$ with $\alpha \geq 0$, then there are no integer solutions to Eq.~\eqref{eq:xppbyppczp} with $xyz \neq 0$.
\end{theorem}

\subsection{General $a,b,c$}

Dahmen and Yazdani~\cite{Dahmen_2012} considered Eq.~\eqref{eq:axppbypmczp} with $c=16$ and
$p \geq 5$, that is,
\begin{equation}
\label{eq:axppbypm24zp}
 ax^p + by^p =16 z^p,
\end{equation}
and their results imply the following theorem.
\begin{theorem}
Let
\begin{equation*}
 (a, b) \in \{(5^2,23^4), (5^8,  37), (5^7, 59^7), (7, 47^7),(11, (5\cdot 17)^2)\}.
\end{equation*}
 Then for odd $p \geq 5$, Eq.~\eqref{eq:axppbypm24zp} has no non-zero integer solutions.
\end{theorem}

In 2002, Kraus~\cite[Proposition 2.3]{Kraus2002} considered Eq.~\eqref{eq:axppbypmczp} with
$a=3$, $b=4$ and $c=5$ and proved the following theorem.

\begin{theorem}
Let $p$ be a prime number congruent to $13$ modulo $24$. Then the equation
\begin{equation}
\label{eq:3xpp4ypp5zp}
  3x^p+4y^p=5z^p
\end{equation}
 has no non-trivial integer solutions.
\end{theorem}

In 2016, Freitas and Kraus~\cite{FREITAS2016751} improved this result and obtained the following.

\begin{theorem}
Let $p \geq 5$ be a prime satisfying $p \equiv 5$ (mod 8) or $p \equiv 19$ (mod 24). Then Eq.~\eqref{eq:3xpp4ypp5zp} has no non-trivial integer solutions.
\end{theorem} 

In the same paper, Freitas and Kraus~\cite{FREITAS2016751} also studied Eq.~\eqref{eq:axppbypmczp}
with $a=3$, $b=8$ and $c=21$ and obtained the following result.
\begin{theorem}
Let $p \geq 5$ be a prime\footnote{The paper~\cite{FREITAS2016751} states the result for
$p >5$, but direct verification shows that it remains correct for $p=5$.} satisfying
$p \equiv 5$ (mod 8) or $p \equiv 23$ (mod 24). Then the equation
\begin{equation*}
 3x^p + 8y^p = 21z^p
\end{equation*}
 has no non-trivial integer solutions.
\end{theorem}

Dieulefait and Soto~\cite{MR4203704} studied Eq.~\eqref{eq:axppbypmczp} under some restrictions on the
prime factors of $abc$, and they proved the following results.\footnote{These results state that
$p$ is sufficiently large. Specifically, $p > G(a,b,c)$ which is an explicit bound dependent on
$a,b,c$.}

\begin{theorem}
Assume that all prime factors of $abc$ are $1$ modulo $3$. Then the equation
\begin{equation}
\label{eq:axppbypm16czp}
  ax^p + by^p = 16cz^p \quad 
\end{equation}
 has no solutions with $xyz \neq 0$ for all sufficiently large $p$.
\end{theorem}

\begin{theorem}
Let $n$ be a positive integer not dividing $14$, $16$ or $18$. Assume that all prime factors of $a,b,c$ are  equal to $\pm 1$ modulo $n$. Then Eq.~\eqref{eq:axppbypm16czp} has no solutions with $xyz \neq 0$ for all sufficiently large $p$.
\end{theorem}

\begin{theorem}
Assume that all prime factors of $abc$ are $1$ modulo $12$. Then for every integer $r\geq 0$, $r\neq 1$, the equation
\begin{equation*}
 ax^p + by^p = 2^rcz^p \quad 
\end{equation*}
 has no solutions with $xyz \neq 0$ for all sufficiently large $p$.
\end{theorem}

\begin{theorem}
\label{th:qoddprime}
 Let $q$ be an odd prime. Assume that $\gcd(a,b,c)=1$, all odd prime factors of $a,b,c$ are equal to $1$ modulo $4q$, and either $abc$ is odd or $4\vert bc$. Then Eq.~\eqref{eq:axppbypmczp} has no solutions with $xyz \neq 0$ for all sufficiently large $p$.
\end{theorem}

In the same work, the authors~\cite{MR4203704} also prove the following result which uses the Legendre
symbol\footnote{For the definition, see~\cite[Sec. 1.4.2]{MR1228206}. } $\left(\frac{\cdot}{\cdot} \right)$.

\begin{theorem}
\label{th:jacobi}
 Let $q \geq 5$ and $l \geq 5$ be primes. Assume one of the following:
\begin{itemize}
  \item[(i)] {$(q,l) \equiv (-5, 5)$} or $(11, -11)$ modulo $24$,
  \item[(ii)] {$q \equiv 11$} modulo $24$, $l \equiv 5$ modulo $24$ and $\left(\frac{q}{l}\right)=-1$, or
  \item[(iii)] {$q \equiv \pm 3$} modulo $8$, $l \equiv -1$ modulo 24, $l \not \equiv -1$ modulo $q$.
\end{itemize}
 Assume that $\gcd(a,b,c)=1$ and $\mathrm{rad}(abc) = ql$. Let $n = 0$ or $n \geq 4$, then the equation
\begin{equation*}
 ax^p + by^p + 2^ncz^p = 0
\end{equation*}
 has no solutions with $xyz \neq 0$ for all sufficiently large $p$.
 Let $r \geq 1$, then the equation
\begin{equation*}
 ax^p + 2^rby^p + 2^rcz^p = 0
\end{equation*}
 has no solutions with $xyz \neq 0$ for all sufficiently large $p$.
\end{theorem}

Powell~\cite{POWELL198434} studied Eq.~\eqref{eq:axppbypmczp} and obtained the following result.  

\begin{theorem}
For any even integer $m$ for which\footnote{Here $\phi(c)$ denotes Euler's totient function.} $3\phi(m) > m$, if $n$ is any positive integer sufficiently large for which $mn + 1 = q$, with $q$ a prime, and $a,b$, and $c$ are integers for which $a \pm b \pm c \neq 0$ and $a \neq \pm b$, $b \neq \pm c$, $a \neq \pm c$, then the equation
\begin{equation}
\label{eq:axnpbynmczn}
  ax^n + by^n = cz^n
\end{equation}
 has no solutions in integers $x,y,z$ with $xyz \neq 0$. 
\end{theorem}

The work~\cite{MR779375} then provides explicit bounds for what it means for $n$ to be
sufficiently large, and obtained the following result.  
\begin{theorem}
\label{th:powell}
 Let $a,b, c$ be non-zero integers such that $a\pm b \pm c \neq 0$,  $a\pm b \neq 0$, $a \pm c \neq 0$, $b \pm c \neq 0$. 
\begin{itemize}
  \item[(a)] Let {$l$} be a prime, $l> | a|  +| b| +| c|$, let $n \geq 1$ be an integer such that $2ln+1=q$ is a prime bigger than $(| a|  +| b| +| c| )^{l-1}$. 
  \item[(b)] Let {$m\geq 2$} be an even integer such that $3\phi(m) > m$,  let $n \geq 1$ be an integer such that $mn +1 = q$ is a prime bigger than $(| a| +| b|  +| c| )^{\phi(m)}$.
\end{itemize}
 If (a) or (b) holds, then Eq.~\eqref{eq:axnpbynmczn} has no solutions in integers $x,y,z$ with $xyz \neq 0$. 
\end{theorem}

Browning and Dietmann~\cite{MR2537095} proved that for any fixed $n>1$ the set of integers
$(a,b,c)$ for which Eq.~\eqref{eq:axnpbynmczn} is solvable in non-zero integers $(x,y,z)$
has density $0$.

\section{Equations of signature $(p,p,r)$}

In this section we discuss the case $p=q$ in \eqref{eq:genfermat}, that is, equation
\begin{equation}
\label{eq:axppbypmczr}
 ax^p+by^p=cz^r,
\end{equation}
where $a,b,c$ are positive integers, $p\neq r$ primes and $p \geq 3$. We discuss the works investigating infinite families of \eqref{eq:axppbypmczr} in which some of $a,b,c$ are fixed whilst $p$ or $r$ varies.

\subsection{The case $a=b=1$}
 
We start with the case $a=b=1$. In this instance Eq.~\eqref{eq:axppbypmczr} reduces to
\begin{equation}
\label{eq:genfermata1b1peq}
 x^p+y^p=cz^r.
\end{equation} 
We remark that for Eq.~\eqref{eq:genfermata1b1peq} Proposition \ref{prop:paircop} is applicable if and only if $c$ is not divisible by a $p$th power of any prime. If this condition is satisfied, we may talk about solutions in pairwise coprime integers $(x, y,z)$. From symmetry, we may also assume that $x>y$.

In 2023, Freitas and Najman~\cite{MR4683842} considered Eq.~\eqref{eq:genfermata1b1peq} where
$p,r > 3$ are primes and $c$ is an odd positive integer, not an $r$th power and
$r \nmid c$. They obtained that for a set of prime exponents $p$ of positive density, the
equation has no non-trivial primitive solutions $(x, y, z)$ such that $2 \mid x + y$ or $r \mid x + y$.

In 2024, Kara, Mocanu and \"Ozman~\cite{kara2024} prove some results conditional on the weak
Frey--Mazur conjecture~\cite[Conjecture 3.9]{kara2024} and the Eichler--Shimura
conjecture~\cite[Conjecture 3.11]{MR1638494,kara2024}. 

\begin{theorem}
Let $p > 5$ be a fixed prime and $c$ be a non-zero integer with $\gcd(c, p) = 1$. Assume the Weak Frey-Mazur Conjecture. Then,
\begin{itemize}
  \item[(i)] If {$p \equiv 3$} modulo $4$, then there exists a (non-effective) constant ${\mathcal{B}}_p$ such that for any prime $r$ with $r > {\mathcal{B}}_p$, Eq.~\eqref{eq:genfermata1b1peq} has no non-trivial primitive integer solutions.
  \item[(ii)] If {$p \equiv 1$} modulo $4$, then there exists a (non-effective) constant ${\mathcal{B}}'_p$ such that for any prime $r$ with $r > {\mathcal{B}}'_p$, Eq.~\eqref{eq:genfermata1b1peq} has no non-trivial primitive integer solutions $(x, y, z)$ where $p \vert z$.
\end{itemize}
\end{theorem}

\begin{theorem}
Let $p > 5$ be a fixed prime such that $p \equiv 1$ modulo $4$ and $c$ be a non-zero integer with $\gcd(c, p) = 1$. Assume the Weak Frey-Mazur Conjecture and the Eichler-Shimura Conjecture. Then there exists a (non-effective) constant ${\mathcal{B}}''_p$ such that for any prime $r$ with $r > {\mathcal{B}}''_p$, Eq.~\eqref{eq:genfermata1b1peq} has no non-trivial primitive integer solutions.
\end{theorem}

\subsubsection{Eq.~\eqref{eq:genfermata1b1peq}
 with fixed $p$}

In this section, we consider Eq.~\eqref{eq:genfermata1b1peq} with $p\geq 3$ fixed and $r$ varies. In this section, we will not assume that $p$ is a prime. We start with the case $p=3$.   

Eq.~\eqref{eq:genfermata1b1peq} with $p=3$, $c=2$ and $r$ even was solved by
Zhang and Bai~\cite{zhang2014}.
\begin{theorem}
\label{th:zhang2014}
 Let $r \geq 2$ be an integer. Then the equation
\begin{equation*}
 x^3+y^3=2z^{2r}
\end{equation*}
 has no integer solutions with $x,y$ coprime other than $(x,y,z)=(1,1,\pm 1)$, $\pm (1,-1,0)$.
\end{theorem}

Eq.~\eqref{eq:genfermata1b1peq} with $p=3$, $c=s^\alpha$ with $s$ prime, $\alpha \geq 1$ and $r\geq 4$, that is,
\begin{equation}
\label{eq:x3y3calphazr}
 x^3+y^3=s^\alpha z^r
\end{equation}
has been studied by several authors. In 2006, Mulholland~\cite{mulholland2006elliptic} obtained the
following result.
\begin{theorem}
\label{th:mulholland2006}
 There is an explicit set $T$ of primes containing a positive proportion of all primes,\footnote{Specifically, $T$ is the set of primes $s$ for which there does not exist an elliptic curve with rational 2-torsion and conductor $2^M 3^{2}s$, $1 \leq M \leq 3$.} such that for every 
 $s \in T$, every integer $\alpha \geq 1$, and every prime $r$ satisfying $r \geq s^{8s}$ and $r \nmid \alpha$, Eq.~\eqref{eq:x3y3calphazr} has no solutions in coprime non-zero integers $x,y,z$.
\end{theorem}

Theorem \ref{th:mulholland2006} has been significantly strengthened by Bennett, Luca and
Mulholland~\cite{bennett2011twisted}. 

\begin{theorem}
\label{th:bennett2011}
 There exists a set ${\mathcal{W}}_3$ of primes, with $\#\{s\leq x: s\in {\mathcal{W}}_3\} \leq K\sqrt{x}\log^2(x)$ for some $K>0$, such that for every prime $s\geq 5$ with $s \not \in {\mathcal{W}}_3$, any $\alpha\geq 1$, and every prime $r$ satisfying $r\geq s^{2s}$, Eq.~\eqref{eq:x3y3calphazr} has no solutions in coprime non-zero integers $x,y,z$. 
\end{theorem}

Because the number of all primes up to $x$ is about $\frac{x}{\log x}$, set ${\mathcal{W}}_3$ contains a
negligible part of all primes, hence Theorem \ref{th:bennett2011} solves Eq.~\eqref{eq:x3y3calphazr}
for almost all primes $s$. There is an algorithm for checking whether any particular prime
belongs to ${\mathcal{W}}_3$. For example, it is known~\cite{bennett2011twisted} that all primes
$s$ equal to $317$ or $1757$ modulo $2040$ do not belong to ${\mathcal{W}}_3$,
hence Theorem \ref{th:bennett2011} is applicable to all such primes.

Theorem \ref{th:bennett2011} is not applicable for primes $s \in {\mathcal{W}}_3$. In 2018, Bennett, Bruni and
Freitas~\cite{MR3830208} proved a more general version of this theorem with a smaller exceptional set.

Let
\begin{equation*}
S_1 = \{q \text{ prime} : q = 2^a 3^b \pm 1, a \in \{2, 3\} \text{ or } a \geq 5, b \geq 0\},
\end{equation*}
\begin{equation*}
S_2 = \{q \text{ prime} : q = | 2^a \pm 3^b| , a \in \{2, 3\} \text{ or } a \geq 5, b \geq 1\},
\end{equation*}
\begin{equation*}
S_3 = \left\{q \text{ prime} : q = \frac{1}{3} (2^a+ 1),a \geq 5 \text{ odd}\right\},
\end{equation*}
\begin{equation*}
S_4 = \{q \text{ prime} : q = d^2+ 2^a 3^b, a \in \{2, 4\} \text{ or } a \geq 8 \text{ even}, b \text{ odd}\},
\end{equation*}
\begin{equation*}
S_5 = \{q \text{ prime} : q = 3d^2+2^a, a \in \{2, 4\} \text{ or } a \geq 8\text{ even}, d \text{ odd}\},
\end{equation*}
\begin{equation*}
S_6 = \left\{q \text{ prime} : q = \frac{1}{4} (d^2+3^b), d \text{ and } b \text{ odd}\right\},
\end{equation*}
\begin{equation*}
S_7 = \left\{q \text{ prime} : q = \frac{1}{4} (3d^2+1), d \text{ odd}\right\},
\end{equation*}
\begin{equation*}
S_8 = \left\{q \text{ prime} : q = \frac{1}{2} (3v^2-1),  u^2-3v^2=-2 \right\},
\end{equation*}
where $a,b,u,v$ and $d$ are integers, and let
\begin{equation}
\label{eq:S0def}
 S_0 = S_1 \cup S_2 \cup S_3 \cup S_4 \cup S_5 \cup S_6 \cup S_7 \cup S_8.
\end{equation}

\begin{theorem}
\label{th:bennett2017}
 For every prime $s\geq 5$ such that $s \not \in S_0$, any $\alpha\geq 1$, and any prime $r$ satisfying $r\geq s^{2s}$, Eq.~\eqref{eq:x3y3calphazr} has no solutions in coprime non-zero integers $x,y,z$.
\end{theorem}

As mentioned in~\cite{MR3830208}, $S_0 \subset {\mathcal{W}}_3$, hence Theorem \ref{th:bennett2017} is a generalization of
Theorem \ref{th:bennett2011}. The smallest examples of primes in ${\mathcal{W}}_3 \setminus S_0$ are 
\fontsize{7.5pt}{10.5pt}\selectfont\begin{equation*}
s = 53, 83, 149, 167, 173, 199, 223, 227, 233, 263, 281, 293, 311, 347, 353, 359, \dots
\end{equation*}\normalsize 
Examples of infinite families outside $S_0$ are primes $s$ equal to $53$ modulo $96$, $120$ or $144$, or primes $s$ equal to $65$ modulo $81$ or $84$. So, Theorem \ref{th:bennett2017} is applicable to all such primes.

Bennett, Bruni and Freitas~\cite{MR3830208} considered the equation
\begin{equation}
\label{eq:x3py3mszr}
 x^3+y^3=sz^r
\end{equation}
where $s$ and $r$ are primes. For small primes $s \in S_0$, they
determined~\cite[Theorem 7.2]{MR3830208} the modularity conditions on $r$ such that for
$r \geq s^{2s}$ the equation has no primitive integer solutions, see 
Table \ref{tb:MR3830208}.

\begin{table}
\begin{center}
  \caption{\label{tb:MR3830208}Modularity conditions on $r$ for small primes $s$, such that Eq.~\eqref{eq:x3py3mszr} has no primitive integer solutions if $r \geq s^{2s}$.}
\begin{tabular}{|c|c||c|c|}
\hline

   $s$ & Condition on $r$ & $s$ & Condition on $r$ \\\hline\hline
   5 & 13,19 or 23 (mod 24) & 47 & $5, 11, 13, 17, 19$ or 23 (mod 24) \\\hline
   11 & $13, 17, 19$ or 23 (mod 24) & 59 & $5, 7, 11, 13, 19$ or 23 (mod 24) \\\hline
   13 & 11 (mod 12) & 67 & $7, 11, 13, 29, 37, 41, 43,59, $ \\ &&& $ 67, 71, 89, 101$ or 103 (mod 120) \\\hline
   17 & 5,17 or 23 (mod 24) & 71 & 5 (mod 6) \\\hline
   23 & 19 or 23 (mod 24) & 73 & 41,71 or 89 (mod 120) \\\hline
   29 & $7, 11, 13, 17, 19$ or 23 (mod 24) & 79 & $5, 7, 11, 13, 19$ or 23 (mod 24) \\\hline
   31 & 5 or 11 (mod 24) & 89 & $13, 17, 19$ or 23 (mod 24) \\\hline
   41 & $5, 7, 11, 17, 19$ or 23 (mod 24) & 97 & 11 (mod 12) 
\\\hline
\end{tabular}
\end{center}
\end{table}
In the case $s=5$, the authors were able to remove the condition $r \geq s^{2s}$ and prove the following theorem.
\begin{theorem}
If $r$ is prime with $r \equiv 13, 19$ or $23$ (mod $24$), then there are no non-trivial primitive integer solutions $(x,y,z)$ to the equation
\begin{equation*}
 x^3+y^3=5z^r.
\end{equation*}
\end{theorem}

Set $S_0$ contains all primes of the form $s=2^a 3^b-1$ for integers $a\geq 5$ and
$b\geq 1$, hence these primes are not covered by Theorem \ref{th:bennett2017}. For such primes, we have
the following result~\cite{MR3830208}, which uses a Jacobi symbol\footnote{For the definition and its
main properties, see~\cite[Section 1.4.2]{MR1228206}.} condition.

\begin{theorem}
\label{th2:bennett2017}
 Suppose that $s = 2^a 3^b-1$ is prime, where $a \geq 5$ and $b \geq 1$ are integers and let $\alpha \geq 1$. If $r > s^{2s}$ is prime and integers $(\alpha,r,a,b)$ do \emph{not} satisfy the condition
\begin{equation}
\label{cd:th2:bennett2017}
  \left( \frac{\alpha}{r} \right) = \left( \frac{4-a}{r} \right) = \left( \frac{-6b}{r} \right),
\end{equation}
 then Eq.~\eqref{eq:x3y3calphazr} has no non-trivial primitive solutions.
\end{theorem}

Another result from~\cite{MR3830208} is the following one.
\begin{theorem}
\label{th3:bennett2017}
 Let $T = S_7 \cup \{q \text{ prime} : q = 3d^2 + 16, d \in \mathbb{Z}\}$. 
 If $s$ is a prime with $s \not \in T$, then, for a positive proportion of primes $r$, there are no solutions to equation
\begin{equation*}
 x^3+y^3=sz^r
\end{equation*} 
 in coprime non-zero integers $x$, $y$ and $z$.
\end{theorem}

We next discuss Eq.~\eqref{eq:genfermata1b1peq} with $p=4$ and $c,r$ primes, that is
\begin{equation}
\label{eq:x4y4czr}
 x^4+y^4 = cz^r.
\end{equation}
If $c$ is not equal to $1$ modulo $8$, then it is easy to show
(see~\cite{MR2139003}) that \eqref{eq:x4y4czr} has no solutions modulo $c$ other than
$(0,0,0)$ and therefore has no primitive integer solutions for any $r$. Hence, it
suffices to consider the case $c \equiv 1(\text{mod}\, 8)$.

We first discuss Eq.~\eqref{eq:x4y4czr} with $c=73,89$ or $113$ and $r >13$.
Dieulefait~\cite{MR2139003} proved the following theorem.
\begin{theorem}
If $r$ is a prime with $r >13$, then Eq.~\eqref{eq:x4y4czr} with $c = 73, 89$ or $113$ does not have any non-trivial primitive integer solutions.
\end{theorem}
The techniques of~\cite{MR2139003} are applicable only to prove that an equation of the form
\eqref{eq:x4y4czr} has no primitive integer solutions. If $c=A^4+B^4$ for some integers
$A,B$, then the equation has an integer solution with $z=1$, so the techniques
of~\cite{MR2139003} are not applicable. For a different reason, the techniques are also not applicable
if $c=(2A)^4 + B^2$ for some integers $A,B$. The values $c = 73, 89, 113$ are the first three values
of $c$ satisfying $c \equiv 1(\text{mod}\, 8)$ and not belonging to these excluded families.

We next discuss Eq.~\eqref{eq:genfermata1b1peq} with $p=5$, that is,
\begin{equation}
\label{eq:x5py5mczr}
 x^5+y^5=cz^r.
\end{equation}

In 2015, Bruni~\cite{bruni2015twisted} proved a version of Theorem \ref{th:bennett2011} with exponent
$5$ instead of $3$.
\begin{theorem}
\label{th:bruni2015}
 There exists a set ${\mathcal{W}}_5$ of primes, with $\#\{s\leq x: s\in {\mathcal{W}}_5\} \leq K\sqrt{x}\log^2(x)$ for some $K>0$, such that for every prime $s\geq 5$ such that $s \not \in {\mathcal{W}}_5$, any $\alpha\geq 1$, and any prime $r$ satisfying $r\geq s^{13s}$, equation
\begin{equation}
\label{eq:x5y5calphazr}
  x^5+y^5=s^\alpha z^r
\end{equation}
 has no solutions in coprime non-zero integers $x,y,z$.
\end{theorem}

Bruni~\cite{bruni2015twisted} also developed methods for solving Eq.~\eqref{eq:x5y5calphazr} for
certain primes $s$ belonging to ${\mathcal{W}}_5$. In fact, he solved Eq.~\eqref{eq:x5y5calphazr} for all primes $s$ outside the specific families listed in Table
7.6 of~\cite{bruni2015twisted}.

For Eq.~\eqref{eq:x5py5mczr} with $c=2$, that is,
\begin{equation}
\label{eq:x5py5mz2}
 x^5+y^5=2z^r,
\end{equation}
we have the following partial result.
\begin{theorem}[\cite{signature2019multi}]
 For all primes $r$, Eq.~\eqref{eq:x5py5mz2} has no non-trivial primitive integer solutions $(x,y,z)$ such that $z$ is divisible by either $2$, $5$, or $r$.
\end{theorem}

Eq.~\eqref{eq:x5py5mczr} with $c=3$, that is,
\begin{equation}
\label{eq:sig55pc3}
 x^5+y^5=3z^r
\end{equation}
has been solved for certain values of $r$
in~\cite{Billerey_2007,billerey2008solvingfermattypeequationsx5y5dzp}, and then for all values of
$r$ by Billerey, Chen, Dieulefait and Freitas~\cite{signature2019multi}.

\begin{theorem}
\label{th:billerey2024}
 Eq.~\eqref{eq:sig55pc3} has no non-trivial primitive integer solutions for any integer $r \geq 2$.
\end{theorem}

We next discuss Eq.~\eqref{eq:x5py5mczr} for other values of $c$.  
In 2013, Dieulefait and Freitas~\cite{44370229} proved the following theorem. 
\begin{theorem}
Let $s$ be an integer divisible only by primes $s_i \not \equiv 1$ (mod 5).
 
\begin{itemize}
  \item[(a)] For any prime {$r>13$} such that $r \equiv 1$ (mod 4), or $r \equiv \pm 1$ (mod 5), Eq.~\eqref{eq:x5py5mczr} with $c=2 s$ has no non-trivial integer solutions with $\gcd(x,y)=1$.
  \item[(b)] For any prime {$r>73$} such that $r \equiv 1$ (mod 4), or $r \equiv \pm 1$ (mod 5), Eq.~\eqref{eq:x5py5mczr} with $c=3 s$ has no non-trivial integer solutions with $\gcd(x,y)=1$.
\end{itemize}
\end{theorem}

In 2007, Billerey~\cite{Billerey_2007} proved the following result. 

\begin{theorem}
Let $c=2^{\alpha} \cdot 5^{\gamma}$ with $2\leq \alpha \leq 4$ and $0 \leq \gamma \leq 4$ and let $r$ be a prime satisfying $r\geq 7$. Then, Eq.~\eqref{eq:x5py5mczr} has no non-trivial coprime integer solutions. 
\end{theorem}

Billerey and Dieulefait~\cite{billerey2008solvingfermattypeequationsx5y5dzp} studied the cases of
$c=7$ and $c=13$ in \eqref{eq:x5py5mczr} and proved the following results.
\begin{theorem}
If either (i) $c=7$ and $r\geq 13$ or (ii) $c=13$ and $r \geq 19$, then Eq.~\eqref{eq:x5py5mczr} does not have any non-trivial primitive integer solutions.
\end{theorem}

In the same paper~\cite{billerey2008solvingfermattypeequationsx5y5dzp}, the authors also proved that Eq.~\eqref{eq:x5py5mczr} with $c=2^{\alpha} 3^{\beta} 5^{\gamma}$, that is,
\begin{equation*}
x^5+y^5= 2^{\alpha} 3^{\beta} 5^{\gamma} z^r
\end{equation*}
with $\alpha \geq 2$, $\beta,\gamma$ arbitrary and $r \geq 13$, has no integer solutions with
$\gcd(x,y)=1$ and $z \neq 0$. Work~\cite{Billerey_2007} studies this equation for primes
$r \geq 7$ and provides sufficient conditions for $r$, $\alpha,\beta$ and $\gamma$
for which the equation has no non-trivial coprime integer solutions.

None of the theorems above include the case when $c$ is divisible by $11$ in
\eqref{eq:x5py5mczr}. Noubissie~\cite{noubissie2020generalized} considered Eq.~\eqref{eq:x5py5mczr}
with $c=2^{\alpha}11$, that is,
\begin{equation}
\label{eq:x5py5m2a11zr}
 x^5+y^5= 2^{\alpha}11 z^r
\end{equation}
and proved the following theorem.
\begin{theorem}
\label{th:Noubissie}
 Let $r >19$ be prime and $\alpha \geq 2$. Then Eq.~\eqref{eq:x5py5m2a11zr} has no non-trivial primitive integer solutions.
\end{theorem}
Noubissie~\cite{noubissie2020generalized} also proved some results in the cases $\alpha=0$ and
$\alpha=1$ which are not included in Theorem \ref{th:Noubissie}. If $\alpha =0$, the resulting equation
$x^5+y^5= 11 z^r$ has no non-trivial primitive solutions with $x\equiv 2(\text{mod}\, 4)$, while if $\alpha=1$ then
the resulting equation $x^5+y^5= 22 z^r$ has no non-trivial primitive solutions with $z$ even.

Let us now discuss Eq.~\eqref{eq:genfermata1b1peq} with $p>5$. In case $p=6$,
Kraus~\cite{kraus2002question} proved the following result.

\begin{theorem}
Let $r\geq 5$ be a prime number. Then equation
\begin{equation*}
 x^6+y^6=2z^r
\end{equation*}
 has no non-trivial primitive integer solutions with $| xyz|  \neq 1$. 
\end{theorem}

For odd $p$, equation
\begin{equation}
\label{eq:kraus2002}
 x^p-y^p=cz^r
\end{equation}
reduces to \eqref{eq:genfermata1b1peq} by the trivial change of variables $y \to -y$. For even
$p$, these equations are not equivalent, and \eqref{eq:kraus2002} must be considered
separately. Kraus~\cite{kraus2002question} studied Eq.~\eqref{eq:kraus2002} with $p$
being an integer in the range $6\leq p \leq 12$. To state the results, we must first introduce the
following notation. 

Given an integer $N \geq 1$, we denote by ${\mathcal{D}}_N$ the set of integers $A \geq 1$ which satisfy both conditions
\begin{itemize}
 \item[(i)] for any prime number {$l$} we have\footnote{For a prime $l$, $\nu_l(A)$ is the $l$-adic valuation of $A$, that is, the highest power of $l$ that divides $A$.} $\nu_l(A) < N$, and,
 \item[(ii)] for any prime divisor {$l$} of $A$ we have $l \not \equiv 1$ modulo $N$.
\end{itemize}

We may now state the theorems proved by Kraus~\cite{kraus2002question}.
\begin{theorem}
\label{th:kraus2002}
\begin{itemize}
  \item[(1)] Suppose we have
\begin{equation*}
  (p,r) \in \{(6,2),(8,2),(8,3),(9,3),(10,2),(12,2) \}.
\end{equation*}
  Then, for any integer $c$ belonging to ${\mathcal{D}}_p$, Eq.~\eqref{eq:kraus2002} has no non-trivial primitive integer solutions.
  \item[(2)] Suppose we have {$p=5$} or $p=7$. Let $c$ be an integer divisible by $p$ and
  belonging to ${\mathcal{D}}_p$. Then, Eq.~\eqref{eq:kraus2002} with $r=2$ has no non-trivial primitive integer solutions.
\end{itemize}
 
\end{theorem}

\begin{theorem}
Let $c \geq 1$ be an integer such that for any prime divisor $l$ of $c$, we have $\nu_l(c) < 12$ and $l \not \equiv 1$ modulo $4$. Then, for any integer $r \geq 3$, the equation
\begin{equation*}
 x^{12}-y^{12}=cz^r
\end{equation*}
 has no non-trivial primitive integer solutions. 
\end{theorem}

We now return to Eq.~\eqref{eq:genfermata1b1peq}, and next discuss it with $p=7$. Let
$s$ be an integer only divisible by primes $s_i \not \equiv 0,1$ modulo $7$. In 2015,
Freitas~\cite{freitas2015recipes} proved that if either (i) $\alpha \geq 2$, $\beta \geq 0$ and
$\gamma \geq 0$, or (ii) $\alpha =1$, $\beta \geq 1$ and $\gamma \geq 0$, or (iii) $\alpha= 0$,
$\beta \geq 0$ and $\gamma \geq 1$, then Eq.~\eqref{eq:genfermata1b1peq} with $p=7$,
$r \geq 17$ and $c=2^{\alpha} 3^{\beta} 5^{\gamma} s$ has no integer solutions where $z$ is divisible by
$7$ and $| xyz|  > 1$. In the same paper, Freitas~\cite{freitas2015recipes} also proved the
following result.
\begin{theorem}
If $\alpha \geq 2$ or $\beta >0$ or $\gamma >0$, then Eq.~\eqref{eq:genfermata1b1peq} with $p=14$,  $r \geq 17$ and $c=2^{\alpha} 3^{\beta} 5^{\gamma} s$, (where $s$ is an integer only divisible by primes $s_i \not \equiv 0,1$ modulo $7$), has no primitive integer solutions with $| xyz|  > 1$.
\end{theorem}

In 2024, the authors of~\cite{billerey2024darmonsprogramgeneralizedfermat} use the Darmon
program~\cite{darmon2000rigid} to solve Eq.~\eqref{eq:genfermata1b1peq} for $(p,c)=(7,3)$.
\begin{theorem}
For any integer $r\geq 2$, equation
\begin{equation*}
 x^7+y^7=3z^r
\end{equation*}
 has no non-trivial primitive integer solutions.
\end{theorem}

The case $(p,r)=(9,2)$ of \eqref{eq:genfermata1b1peq} has been studied by Cel~\cite{MR995897}, 
who proved the following result.

\begin{theorem}
For any $n \geq 0$, the equation
\begin{equation*}
 x^9+y^9=2^nz^2
\end{equation*}
 has no non-trivial primitive integer solutions with $z \neq \pm 1$.
\end{theorem}

We next discuss the case $p=13$ in \eqref{eq:genfermata1b1peq}. In 2014, Siksek and
Frietas~\cite{freitas2014criteriairreducibilitymodp}, improving an earlier result of Dieulefait and
Freitas~\cite{FreitasDieulefait2013}, proved the following theorem.

\begin{theorem}
Let $s = 3, 5, 7$ or $11$ and $\gamma$ be an integer divisible only by primes $\gamma_i \not \equiv 1$ (mod 13). Let also
\begin{equation}
\label{def:R}
  R \coloneq  \{2, 3, 5, 7, 11, 13, 19, 23, 29,  71\}.
\end{equation}
 If $r$ is a prime not belonging to $R$, then:
\begin{itemize}
  \item[(i)] Eq.~\eqref{eq:genfermata1b1peq} with {$p=13$} and $c=s \gamma$ has no primitive integer solutions $(x,y,z)$ such that $13 \nmid z$ and $| xyz|  > 1$.
  \item[(ii)] Eq.~\eqref{eq:genfermata1b1peq} with {$p=26$} and $c=10 \gamma$ has no primitive integer solutions with $| xyz|  > 1$.
\end{itemize}
\end{theorem}

Billerey, Chen, Dieulefait, and Freitas~\cite{signature2019multi} proved that
Eq.~\eqref{eq:genfermata1b1peq} with $c=3$, $p=13$ and $r \geq 2$ where
$r$ is not a power of $7$ has no non-trivial primitive integer solutions. In a
later paper, Billerey, Chen, Demb\'el\'e, Dieulefait, and Freitas~\cite{10.5565/PUBLMAT6722309}
extended this result to all integers $r \geq 2$.

\begin{theorem}
Equation
\begin{equation*}
 x^{13}+y^{13}=3z^r
\end{equation*}
 has no non-trivial primitive integer solutions for all integers $r \geq 2$.
\end{theorem}

\subsubsection{Eq.~\eqref{eq:genfermata1b1peq}
 with fixed $r$}

We now consider Eq.~\eqref{eq:genfermata1b1peq} with $r$ fixed. We start with the case $r=2$, that is, equation
\begin{equation}
\label{eq:xppypmcz2}
 x^p+y^p=cz^2.
\end{equation}  

In 1977, Terjanian~\cite{00efa9d9-bb1b-3dd3-b2cb-1761e27ea5c8} proved for an odd prime $p$
and integer $c \geq 2$ with no square factor or a prime divisor of the form $2kp+1$, if
$(x,y,z)$ is a non-trivial primitive integer solution to \eqref{eq:xppypmcz2} with even
$z$, then $z$ must be divisible by $p$.

In 2011, S\"oderlund and Ald\'en~\cite{soderlund2011diophantine} studied Eq.~\eqref{eq:xppypmcz2} with $p=5$ and proved the following result.

\begin{theorem}
\label{th:soderlund2011}
 If $\alpha = 0$ or $\alpha \geq 2$, $\beta \geq 0$, $\gamma \geq 0$, $s$ is a prime equal to $3$ or $7$ modulo $10$, then the equation
\begin{equation}
\label{eq:x5py5mcz2}
  x^5+y^5=cz^2
\end{equation}
 with $c = 2^\alpha 5^\beta s^\gamma$ has no non-trivial primitive integer solutions.
\end{theorem}
The examples of values of $c$ covered by Theorem \ref{th:soderlund2011} are
\begin{equation*}
c = 1, 3, 5, 7, 8, 13, 15, 17, 23, 24, 35, 37, 40, 43, 47, 53, 56, 65, 67, 73, \dots
\end{equation*}
Theorem \ref{th:soderlund2011} does not include the case $\alpha=1$. In the follow-up
work~\cite{soderlund2013note}, the same authors included this case as well, and thus solved Eq.~\eqref{eq:x5py5mcz2} for $c=10$ and $c=2s$  where $s$ is a prime equal to
$3$ or $7$ modulo $10$.

Bennett and Skinner~\cite{Bennett_Skinner_2004} considered Eq.~\eqref{eq:xppypmcz2} with
$p \geq 4$, and they obtained the following result. 
\begin{theorem}
\label{th:BennettSkinner2004}
 If $p \geq 4$ and $c \in \{2, 3, 5, 6, 10, 11, 13, 17 \}$ then Eq.~\eqref{eq:xppypmcz2} has no solutions in non-zero pairwise coprime integers $(x, y,z)$ with $x > y$, unless $(p,c) = (4, 17)$ or
\begin{equation*}
 (p,c, x, y,z) \in \{(5, 2, 3, -1, \pm 11),(5, 11, 3, 2, \pm 5),(4, 2, 1, -1, \pm 1)\}.
\end{equation*}
\end{theorem}
Bennett and Skinner~\cite{Bennett_Skinner_2004} state that with further computation, the results can
be extended to cover the additional cases $c=14,15$ and $19$. Some cases of
Theorem \ref{th:BennettSkinner2004} rely on results proven by Bruin in~\cite{Bruin2006}.

Eq.~\eqref{eq:xppypmcz2} has also been studied by Bennett and
Mulholland~\cite{bennett2006diophantine}, who proved that if $p$ and $c \geq 5$ are
primes with $c \neq 7$ and $p > c^{12c^2}$, then  \eqref{eq:xppypmcz2} has no primitive solutions
with $| xy|  > 1$ and $z$ even. They also considered the case  $c=2s$ where
$s \geq 5$ is prime, and proved the following theorem.
\begin{theorem}
If $p$ and $s \geq 5$ are primes such that $p > s^{132s^2}$, then the equation
\begin{equation*}
 x^p + y^p = 2sz^2
\end{equation*}
 has no primitive integer solutions with $| xy|  > 1$.
\end{theorem}

Bennett and Mulholland~\cite{mulholland2006elliptic} studied Eq.~\eqref{eq:xppypmcz2} with
coefficient $c=2^{\alpha}s$ where $s$ is prime and $\alpha \geq 1$, and proved the following
theorem.
\begin{theorem}
Let $s \neq 7$ be prime and $\alpha \geq 1$. If prime $p$ satisfies $p > s^{27s^2}$, then the equation
\begin{equation}
\label{eq:xppypm2asz2}
  x^p + y^p = 2^{\alpha}sz^2
\end{equation}
 has no primitive integer solutions.
\end{theorem}

D\c abrowski~\cite{MR2737959} proved a more general version of this theorem, which does not restrict
$s$ to being a prime.

\begin{theorem}
Let $s$ be an odd square-free positive integer with $\gcd(s,21)=1$ and $\alpha \geq 1$. If prime $p$ satisfies $p >s^{132s^2}$, then Eq.~\eqref{eq:xppypm2asz2} has no primitive integer solutions. 
\end{theorem}

Ivorra~\cite{MR2310336} studied Eq.~\eqref{eq:xppypmcz2} for $p \in \{11,13,17 \}$ and obtained the
following theorem. 

\begin{theorem}
Let $p \in \{11,13,17 \}$ and let $c \geq 3$ be an integer without square factors and for any prime divisor $l$ of $c$ we have $l \not \equiv 1$ modulo $p$. Then Eq.~\eqref{eq:xppypmcz2} has no non-trivial primitive integer solutions.
\end{theorem}

Bennett, Vatsal and Yazdani~\cite{bennett2004ternary} studied Eq.~\eqref{eq:genfermata1b1peq}
with $r=3$, and proved the following theorems.
\begin{theorem}
\label{th2:bennett2004}
 If $c$ and $p$ are integers with $2 \leq c \leq 5$ and $p \geq 4$, then the equation
\begin{equation}
\label{eq:xpypcz3}
  x^p + y^p = cz^3
\end{equation}
 has no solutions in coprime non-zero integers $x$, $y$ and $z$ with $| xy|  > 1$.
\end{theorem} 

\begin{theorem}
\label{th:bennett2004}
 If $s$ and $p$ are primes such that $p > s^{4s^2}$, and $\alpha$ is a non-negative integer, then the equation
\begin{equation*}
 x^p + y^p = s^\alpha z^3
\end{equation*}
 has no solutions in coprime integers $x$, $y$ and $z$ with $| xy|  > 1$.
\end{theorem}

Krawci\'ow~\cite{KarolinaKrawciow2011} studied Eq.~\eqref{eq:xpypcz3} where $c$ is
divisible by $3$ and obtained the following result.

\begin{theorem}
\label{th:Krawciow2011}
 Let $p$ be a prime number, and let $c$ be a non-zero cube-free integer which is divisible by $3$. Let $C=\mathrm{rad}(c)$, then if $p > C^{10C^2}$, Eq.~\eqref{eq:xpypcz3} has no non-trivial primitive integer solutions.
\end{theorem}

\subsection{The case $a=c=1$}

In this case, Eq.~\eqref{eq:axppbypmczr} is of the form
\begin{equation}
\label{eq:genfermata1c1ppr}
 x^p+by^p=z^r.
\end{equation}

Eq.~\eqref{eq:genfermata1c1ppr} with $r=2$, $b=2^{\alpha}$ for $\alpha \geq 2$, and $p \geq 7$
prime was studied by Cohen~\cite[Section 15.3]{MR2312338} who obtained the following result.
\begin{theorem}
\label{th:xpp2aypmz2}
 The only non-zero pairwise coprime solutions to
\begin{equation}
\label{eq:xpp2aypmz2}
  x^p+ 2^{\alpha}y^p = z^2
\end{equation}
 for $\alpha \geq 2$ and $p\geq 7$ prime are for $\alpha=3$, for which $(x,y,z)=(1,1,\pm 3)$ is a solution for all $p$.
\end{theorem}

Theorem \ref{th:xpp2aypmz2} describes pairwise coprime solutions to \eqref{eq:xpp2aypmz2} under the stated
conditions, but this equation can have primitive solutions that are not pairwise coprime, e.g.~$(x, y,z) = (2, 1,3 \cdot 2^{(p-3)/2})$ for $\alpha = p - 3$, see~\cite{ivorra2003equations}. In 2004, Ivorra~\cite[Chapter
1]{ivorra2004equations} proves that there are no other such examples, which completely solves Eq.~\eqref{eq:xpp2aypmz2} for $\alpha \geq 2$ and $p\geq 7$. Siksek also independently proved the same
result in~\cite{MR2142239}.

Theorem \ref{th:xpp2aypmz2} excludes the case where $p=5$. In this case,
Ivorra~\cite{ivorra2003equations} proved the following result. 

\begin{theorem}
If $1 \leq \alpha \leq 4$ then for Eq.~\eqref{eq:xpp2aypmz2} with $p=5$, that is,
\begin{equation*}
 x^5+ 2^{\alpha}y^5 = z^2,
\end{equation*}
 the non-trivial primitive integer solutions $(x,y,z)$ are given by
\begin{equation*}
 (\alpha,x,y,z)=(1,-1,1,\pm 1),(2,2,1,\pm 6),(3,1,1,\pm 3),(4,2,-1,\pm 4).
\end{equation*}
\end{theorem}

Bennett, Vatsal and Yazdani~\cite{bennett2004ternary} considered Eq.~\eqref{eq:genfermata1c1ppr} with
$r=3$ and $b$ being a prime power, that is,
\begin{equation}
\label{eq:bennett2004ternary}
 x^p + s^{\alpha}y^p = z^3,
\end{equation}
and they obtained the following results.

\begin{theorem}
\label{th:bennett2004ternary}
 If $s$ and $p$ are primes such that $s \neq s_1^3 \pm 3^{s_2}$ for any integers $s_1$ and $s_2$ with $s_2 \neq 1$, $p > s^{2s}$ and $\alpha \geq 0$, then Eq.~\eqref{eq:bennett2004ternary}
 has no solutions in coprime integers $x$, $y$ and $z$ with $| xy|  > 1$.
\end{theorem}

\begin{theorem}
\label{th6:bennett2004}
 Let $\alpha \geq 1$ be an integer,
\begin{equation*}
 s \in \{5, 11, 13, 23, 29, 31, 41, 43, 47, 53, 59, 61, 67, 71, 79, 83, 97\}
\end{equation*}
 and let $p\geq 11$ be a prime not dividing $s^2-1$. Assume also that
\begin{equation*}
\begin{aligned}
 (s, p) \not\in \{  & (13, 19),(29, 11),(43, 13),(47, 13),(59, 11), \\ & (61, 61),(67, 73),(79, 97),(97, 13),(97, 79)\}.
\end{aligned}
\end{equation*} 
 Then Eq.~\eqref{eq:bennett2004ternary} has no solutions in coprime integers $x$, $y$ and $z$ with $| xy|  > 1$. 
\end{theorem}

So far in this section we considered the cases $r=2$ and $r=3$. For the next case
$r=5$, we have the following result~\cite[Theorem 2.8]{soderlund2019some}.

\begin{theorem}
If $b \in \{3, 11, 19, 43, 67, 163\}$, then the equation
\begin{equation*}
 x^4 + b y^4 = z^5
\end{equation*}
 has no non-trivial primitive solutions.
\end{theorem}

Azon \cite{azon2025} considered Eq.~\eqref{eq:genfermata1c1ppr} where $r$ varies. Specifically, he considered Eq.~\eqref{eq:genfermata1c1ppr} with $b=7$, $p=5$ and $r>41$ a prime and he proved that the equation has no non-trivial primitive integer solutions such that $10|z$.

\subsection{The case $a=1$}

In this instance, Eq.~\eqref{eq:axppbypmczr} reduces to
\begin{equation}
\label{eq:xppbyppczr}
 x^p+by^p=cz^r,
\end{equation}
with $\min\{b,c\} > 1$.

We first consider the case $r=2$, which reduces Eq.~\eqref{eq:xppbyppczr} to
\begin{equation}
\label{eq:xppbyppcz2}
 x^p+by^p=cz^2.
\end{equation}

Ivorra~\cite{ivorra2003equations} studies Eq.~\eqref{eq:xppbyppcz2} with $b= 2^{\beta}$ for
$0 \leq \beta \leq p-1$, $c=2$ and $p \geq 7$ prime, that is,
\begin{equation}
\label{eq:ivorra2003}
 x^p+2^{\beta}y^p=2z^2,
\end{equation}
and obtained the following result.

\begin{theorem}
For prime $p \geq 7$ and integer $1 \leq \beta \leq p-1$, Eq.~\eqref{eq:ivorra2003} has no non-trivial primitive integer solutions.
\end{theorem}
Let prime $p \geq 7$ and $\alpha \geq 2$. The task of finding all primitive integer solutions to Eq.~\eqref{eq:xppbyppcz2} with $b=2^{\alpha}$ and $c=3$, that is,
\begin{equation*}
x^p+2^{\alpha} y^p=3z^2,
\end{equation*}
is left as an exercise to the reader in~\cite[Section 15.8]{MR2312338}. The authors
of~\cite{MR3222262} prove that the equation has no solutions in coprime non-zero integers
$x,y,z$ with $| xy|  >1$.

Bennett and Skinner~\cite{Bennett_Skinner_2004} obtained the following results concerning Eq.~\eqref{eq:xppbyppcz2}.
\begin{theorem}
Suppose that $p \geq 7$ is prime. Let $b=2^{\alpha}$. If
\begin{equation}
\label{eq:BennettSkinner2004}
  (c, \alpha_0) \in \{(3, 2),(5, 6),(7, 4),(11, 2),(13, 2),(15, 6),(17, 6)\},
\end{equation}
 $\alpha \geq \alpha_0$, $p > c$ and $(c, \alpha, p) \neq (11, 3, 13)$, 
 then Eq.~\eqref{eq:xppbyppcz2}
   has no solutions in non-zero pairwise coprime integers $(x,y,z)$ with $xy \neq \pm 1$.
\end{theorem}

\begin{theorem}
Suppose that $p \geq 11$ is prime and $\alpha$ is a non-negative integer. If $b \in \{5^{\alpha},11^{\alpha},13^{\alpha}\}$ and $p$ is not a divisor of $b$, then Eq.~\eqref{eq:xppbyppcz2} with $c=2$ has no solution in non-zero pairwise coprime integers $(x, y,z)$.
\end{theorem}

Eq.~\eqref{eq:xppbyppcz2} with $b$ a power of $2$ and $c$ a prime was
studied by Zhang~\cite{zhang2012}, who proved the following theorem.

\begin{theorem}
\label{th:zhang2012}
 Let $s,p$ be odd primes. For any integers $k,d,s\neq (2^k \pm 1)/d^2$ and $\alpha \geq 0$, $\alpha \neq 1$. If $p >s^{8s^2}$, then the equation
\begin{equation*}
 x^p+2^\alpha y^p=sz^2
\end{equation*}
 has no integer solutions in pairwise coprime integers $x,y,z$ with $xyz \neq 0$.
\end{theorem}

From Theorem \ref{th:zhang2012}, Zhang~\cite{zhang2012} deduces the following.

\begin{corollary}
Let $s,p$ be odd primes. For any integers $\alpha \geq 2$ and $s>5$ satisfying $s \equiv \pm 3$ modulo $8$. If $p >s^{8s^2}$, then the equation
\begin{equation*}
 x^p+2^\alpha y^p=sz^2
\end{equation*}
 has no integer solutions in pairwise coprime integers $x,y,z$ with $xyz \neq 0$.
\end{corollary}

We next consider the case $r=3$ in \eqref{eq:xppbyppczr}, that is, equation
\begin{equation}
\label{eq:xppbyppcz3}
 x^p+by^p=cz^3.
\end{equation}
Bennett et al.~\cite{bennett2004ternary} studied \eqref{eq:xppbyppcz3} with $b=s^{\alpha}$ and
$c=3^{\beta}$ for prime $s$ and positive integers $\alpha,\beta$, that is,
\begin{equation}
\label{eq3:bennett2004}
 x^p + s^{\alpha}y^p = 3^{\beta}z^3,
\end{equation}
and they proved the following theorems.

\begin{theorem}
\label{th3:bennett2004}
 If $s$ and $p$ are primes such that $s \not\in \{5, 3t^3 \pm1, 9t^3 \pm 1 : t \in \mathbb{N}\}$, $p>s^{28s}$ and $\alpha, \beta$ are positive integers with $3 \nmid \beta$, then Eq.~\eqref{eq3:bennett2004} has no solutions in coprime integers $x$, $y$ and $z$ with $| xy|  > 1$.
\end{theorem}

\begin{theorem}
\label{th7:bennett2004}
 If $p \geq 7$ is prime, $s \in \{7, 11, 13\}$ with
\begin{equation*}
 (s,p) \not \in \{ (7,13),(13,7)\}
\end{equation*}
 and $\alpha, \beta$ are positive integers with $3 \nmid \beta$, then Eq.~\eqref{eq3:bennett2004} has no solutions in coprime integers $x$, $y$ and $z$ with $| xy|  > 1$.
\end{theorem}

Bennett et al.~\cite{bennett2004ternary} also studied equation
\begin{equation}
\label{eq4:bennett2004}
 x^p + 3^{\alpha}y^p = s^{\beta}z^3
\end{equation}
and proved the following theorem. 

\begin{theorem}
\label{th4:bennett2004}
 If $s \in \{2, 3, 5, 7, 11, 13, 15, 17, 19\}$, $p$ is a prime satisfying $p > \max\{s, 4\}$,
\begin{equation*}
 (s,p) \not \in \{(7,11), (11,13)\} \quad \text{and} \quad (\alpha,s) \not \in \{(1,t) : t =2, \text{ or } t\geq 11\},
\end{equation*} 
 and $\alpha$ and $\beta$ are non-negative integers, then Eq.~\eqref{eq4:bennett2004} 
 has no solutions in coprime integers $x$, $y$ and $z$ with $| xy|  > 1$, unless $(| x| , | y| , \alpha, p, | s^{\beta} z^3| ) =(2, 1, 1, 7, 125)$.
\end{theorem}

Krawci\'ow~\cite{KarolinaKrawciow2011} studied Eq.~\eqref{eq:xppbyppcz3} with $b=s^{\alpha}$,
that is,
\begin{equation}
\label{eq:xppsayppcz3}
 x^p+ s^{\alpha} y^p=cz^3,
\end{equation}
and solved it under some conditions on the parameters.  
To state the conditions, we must first introduce some notation. Let $s_1,\ldots , s_k$ be odd primes, such that $3 \in \{s_1,\ldots , s_k\}$. Then let $P(s_1,\ldots , s_k)$ denote the set of primes $x$ satisfying any of the $2^{k+1}$  equations
\begin{equation*}
x^\alpha - {s_1}^{\alpha_1} \cdots {s_k}^{\alpha_k} y^3 = \pm 1, \quad (1\leq \alpha_i \leq 2, \,\, i = 1,\ldots , k), 
\end{equation*}
for some integers $\alpha>1$ and $y>1$. Theorem \ref{th:Krawciow2011} implies that the set
$P(s_1,\ldots , s_k)$ is always finite. For example, $P(3) = \{2,5\}$, see~\cite{KarolinaKrawciow2011}.

\begin{theorem}
\label{th:KarolinaKrawciow2011}
 Let $c = \prod_{i=1}^{k} s_i^{\gamma_i}$ be a prime factorization of a positive cube-free odd integer $c$ which is divisible by $3$, $\alpha$ be a positive integer, $p$ is a prime, and let $C=\mathrm{rad}(c)$. If $s$ is a prime such that $s \not \in P(s_1, \dots, s_k)$ and $s \neq \prod_{i=1}^{k} s_i^{\alpha_i}t^3 \pm 1$ $(1 \leq \alpha_i \leq 2, i = 1,\dots, k)$ for any integer $t$, and if $p > (sC)^{10sC^2}$, then Eq.~\eqref{eq:xppsayppcz3} 
 has no non-trivial primitive integer solutions.
\end{theorem}

\subsection{The case $c=1$}

In this instance, Eq.~\eqref{eq:axppbypmczr} reduces to
\begin{equation}
\label{eq:eq:axppbypmzr}
 ax^p+by^p=z^r
\end{equation}
for positive integers $a,b$.

\subsubsection{Eq.~\eqref{eq:eq:axppbypmzr}
 with fixed $p$}
Eq.~\eqref{eq:eq:axppbypmzr} with $p=3$ was considered by Bennett and Dahmen which led to the
following partial result~\cite[Theorem 13.2]{bennett2013klein}.

\begin{theorem}
\label{th:bennett2013klein}
 Let $a$ and $b$ be odd integers such that $ax^3 + by^3$ is irreducible in $\mathbb{Q}[x, y]$, and suppose that all solutions to the equation
\begin{equation}
\label{eq:genfermatc1prod}
  ax^3 + by^3 = \prod_{s \vert  6ab} s^{\alpha_s}
\end{equation}
 in coprime non-zero integers $x$ and $y$, and non-negative integers $\alpha_s$, satisfy $\alpha_2 \in \{1, 4\}$. Then there exists an effectively computable constant $r_0$ such that the equation
\begin{equation}
\label{eq:genfermatc1z}
  ax^3 + by^3 = z^r
\end{equation}
 has no solution in non-zero integers $x$, $y$ and $z$ (with $\gcd(x, y) = 1$) and prime $r \equiv 1$ (mod 3) with $r \geq r_0$. If all solutions to \eqref{eq:genfermatc1prod} with $x$ and $y$ coprime integers have $\alpha_2 \leq 4$, then there exists an effectively computable constant $r_1$ such that Eq.~\eqref{eq:genfermatc1z} has no solutions in odd integers $x$, $y$ and prime $r \geq r_1$.
\end{theorem}

Examples of pairs $(a,b)$ satisfying the conditions of Theorem \ref{th:bennett2013klein} are $(1,57)$, $(1,83)$ and $(3, 19)$.

\subsubsection{Eq.~\eqref{eq:eq:axppbypmzr}
 with fixed $r$}

In 2004, Bennett and Skinner~\cite{Bennett_Skinner_2004} considered Eq.~\eqref{eq:eq:axppbypmzr}
with $r=2$, and $p \geq 11$ prime, that is,
\begin{equation}
\label{eq:bennetAB}
 ax^p + by^p = z^2,
\end{equation}
and they obtained the following results.

\begin{table}
\begin{center}
   \caption{\label{tb:ABvalsodd}Values for $a,b,p, \alpha$ not covered by Theorem \ref{th:Bennettabc}.}
\begin{tabular}{|c|c|c||c|c|c|}
\hline

    $ab$ & $p$ & $\alpha$ &$ab$ & $p$ & $\alpha$ \\\hline\hline
    $2^{\alpha} 19^{\beta}$ & 11 & 1 & $2^{\alpha} 61^{\beta}$ & 13 & 0,3 \\\hline
    $2^{\alpha} 43^{\beta}$ & 11 & $0,3,\geq 7$ & $2^{\alpha} 61^{\beta}$ & 31 & $\geq 7$ \\\hline
    $2^{\alpha} 53^{\beta}$ & 11 & $2,4,5$ & $2^{\alpha} 67^{\beta}$ & 11 & $0,3,6$ \\\hline
    $2^{\alpha} 53^{\beta}$ & 17 & 0,3 & $2^{\alpha} 67^{\beta}$ & 13 & 0,3 \\\hline
    $2^{\alpha} 59^{\beta}$ & 11 & $\geq 7$ & $2^{\alpha} 67^{\beta}$ & 17 & $0,3,\geq 7$ \\\hline
    $2^{\alpha} 59^{\beta}$ & 29 & 6 &&& 
\\\hline
\end{tabular}
\end{center}
\end{table}

\begin{theorem}
\label{th:Bennettabc}
 Suppose that $p \geq 11$ is prime, $a,b$ are coprime integers satisfying $p \nmid ab$, $\alpha \geq 0$ and $\beta \geq 1$ are integers, and the values $ab, p, \alpha, \beta$ are not those given in 
Table \ref{tb:ABvalsodd}. If
\begin{equation*}
 ab \in  \{2^{\alpha}11^{\beta}, 2^{\alpha}13^{\beta},2^{\alpha}19^{\beta},2^{\alpha}29^{\beta},2^{\alpha}43^{\beta},2^{\alpha}53^{\beta}, 2^{\alpha}59^{\beta}, 2^{\alpha}61^{\beta},2^{\alpha}67^{\beta} \}
\end{equation*}
 for $\alpha =0$ or $\alpha \geq 3$, or if
\begin{equation*}
 ab \in \{2 \cdot 19^{\beta}, 4 \cdot 11^{\beta}, 4 \cdot 19^{\beta}, 4 \cdot 43^{\beta}, 4 \cdot 59^{\beta}, 4\cdot 61^{\beta}, 4\cdot  67^{\beta}\}
\end{equation*}
 then Eq.~\eqref{eq:bennetAB} has no solutions in non-zero pairwise coprime integers $(x, y,z)$.
\end{theorem}

\begin{theorem}
\label{th:Bennettabceven}
 Suppose that $p \geq 11$ is prime, $a$ and $b$ are coprime integers satisfying $p \nmid ab$,  $\alpha,\beta,\gamma$ are
 non-negative integers with $\alpha \geq 6$ and $ab = 2^{\alpha}s^{\beta}t^{\gamma}$ where
\begin{equation*}
\begin{aligned}
  (s,t)  \in \{ & (3, 31) (\beta \geq 1), (5, 11) (\alpha \geq 7), (5, 19),(5, 23) (\beta \geq 1), \\ &
  (7, 19) (\gamma \geq 1),   (11, 13),  (11, 23) (\beta \geq 1), (11, 29),(11, 31) (\beta \geq 1), \\ &
  (13, 31) (\beta \geq 1),  (19, 23) (\beta \geq 1), (19, 29),(29, 31) (\beta \geq 1)\},
\end{aligned}
\end{equation*}
 and the values $ab$, $p$ and $\alpha$ are not those given in 
Table \ref{tb:ABvalseven}.
 Then Eq.~\eqref{eq:bennetAB} has no solutions in non-zero pairwise coprime integers $(x, y, z)$.
\end{theorem}

In 2014, Zhang and Luo~\cite{MR3222262} obtained the following result in the case $ab = 2^{\alpha}3^{\beta}$
of Eq.~\eqref{eq:bennetAB}.

\begin{table}
\begin{center}
   \caption{\label{tb:ABvalseven}Assume that $\beta$ and $\gamma$ are positive integers. Values for $a,b,p, \alpha$ not covered by Theorem \ref{th:Bennettabceven}.}
\begin{tabular}{|c|c|c||c|c|c|}
\hline

    $ab$ & $p$ & $\alpha$ &$ab$ & $p$ & $\alpha$ \\\hline\hline 
    $2^{\alpha} 5^{\beta} 23^{\gamma}$ & 11 & $\geq 7$ & $2^{\alpha} 11^{\beta} 29^{\gamma}$ & 13 & $\geq 7$ \\\hline
    $2^{\alpha} 7^{\beta} 19^{\gamma}$ & 11 & $\geq 7$ & $2^{\alpha} 19^{\beta} 23^{\gamma}$ & 11 & $\geq 7$ \\\hline
    $2^{\alpha} 11^{\beta} 23^{\gamma}$ & 11 & $\geq 6$ & $2^{\alpha} 19^{\beta} 29^{\gamma}$ & 11 & $\geq 7$ 
\\\hline
\end{tabular}
\end{center}
\end{table}

\begin{theorem}
Suppose that $a$ and $b$ are coprime integers with $ab = 2^{\alpha}3^{\beta}$ where $\alpha \geq 6$ and $\beta \geq 0$, and $p \geq 7$ is prime. Then Eq.~\eqref{eq:bennetAB} has no solutions in non-zero pairwise coprime integers $(x, y, z)$.
\end{theorem}

\subsection{The general $a,b,c$}

In this case, we are considering Eq.~\eqref{eq:axppbypmczr} with general positive coefficients $a,b,c$. In this section, we will sometimes consider equations with $p$ or $r$ composite. If $p$ is even, we may also have $b<0$.

We start with the case $r=2$, that is, equation
\begin{equation}
\label{eq:axppbypmcz2}
 ax^p+by^p=cz^2.
\end{equation} 

Noubissie and Togb\'e~\cite{armandphdthesis,Noubissie21} studied Eq.~\eqref{eq:axppbypmcz2} with
$a=5^{\alpha}$, $\alpha \geq 1$, $b=64$, and $c=3$, and obtained the following result.

\begin{theorem}
Suppose that $p \geq 7$ is a prime number and $\alpha \geq 1$. Then the
 equation
\begin{equation*}
 5^{\alpha}x^p + 64y^p = 3z^2
\end{equation*}
 has no non-zero primitive integer solutions with $x y \equiv 1$ (mod 2). 
\end{theorem}

In 2006, Ivorra and Kraus~\cite{ivorra2006quelques} considered Eq.~\eqref{eq:axppbypmcz2} with
\begin{equation}
\label{abc_in_axppbypmcz2}
 (a,b,c)\in \{(2^n,l^m,1),(2^n l^m,1,1),(1,l^m,2)\}.
\end{equation}
They provide some restrictions on $(l,n)$ and for $m \geq 1$, such that for $p$ larger than an explicit bound which depends on $l$, $m$ and $n$, Eq.~\eqref{eq:axppbypmcz2} has no non-trivial primitive integer solutions.

Ivorra and Kraus~\cite{ivorra2006quelques} also solved Eq.~\eqref{eq:axppbypmcz2} for certain
values of $(a,b,c)$ not of the form \eqref{abc_in_axppbypmcz2}. Specifically, they obtain the
following results.

\begin{theorem}
\begin{itemize}
  \item[(a)] For {$p \geq 7$} and $(a,b,c)=(4,1,3)$, Eq.~\eqref{eq:axppbypmcz2} only has the non-trivial primitive integer solutions $(x,y,z)=(1,-1,\pm 1)$.
  \item[(b)] For {$p \geq 11$} and $(a,b,c)=(64,1,7)$, Eq.~\eqref{eq:axppbypmcz2}  only has the non-trivial primitive integer solutions $(x,y,z)=(1,-1,\pm 3)$.
\end{itemize}
\end{theorem} 

In 2022, Cha\l upka et al.~\cite{CHALUPKA2022153} studied Eq.~\eqref{eq:axppbypmcz2} with
$(a,b,c,p)=(3^{31},-4,7,34)$ and proved the following theorem.
\begin{theorem}
The equation
\begin{equation*}
 3^{31}x^{34}-4y^{34} = 7z^2
\end{equation*}
 has no solutions in coprime odd integers.
\end{theorem}

We next discuss Eq.~\eqref{eq:axppbypmczr} with $r=3$, that is,
\begin{equation}
\label{eq:axppbypmcz3}
 ax^p+by^p=cz^3.
\end{equation} 
Noubissie and Togb\'e~\cite{armandphdthesis,Noubissie21} considered Eq.~\eqref{eq:axppbypmcz3}
with $a=2^{\alpha}$, $\alpha \geq 1$, $b=27$, $c=s^{\beta}$, $\beta \geq 1$, and $s \in \{7,13\}$, and
proved the following results. 
\begin{theorem}
\label{th:noubissie1}
 Suppose that $p \geq 11$ is a prime number and $\alpha \geq 1$, $\beta \geq 1$ and $\beta \equiv 1$ (mod 3). Then the
 equation
\begin{equation*}
 2^{\alpha}x^p + 27y^p = 7^{\beta}z^3
\end{equation*}
 has no non-zero primitive integer solutions. 
\end{theorem}

\begin{theorem}
\label{th:noubissie2}
 Suppose that $p > 13$ is a prime number, $\alpha \geq 1$, $\beta \geq 1$ and $\beta \equiv 1$ (mod 3). Then the
 equation
\begin{equation*}
 2^{\alpha}x^p + 27y^p = 13^{\beta} z^3
\end{equation*}
 has no non-zero primitive integer solutions. 
\end{theorem}

We next discuss Eq.~\eqref{eq:axppbypmczr} with $r=5$, that is,
\begin{equation}
\label{eq:axppbypmcz5}
 ax^p+by^p=cz^5.
\end{equation} 
Azon \cite{azon2025} considered Eq.~\eqref{eq:axppbypmcz5} with $a=7$, $b=1$, $c=3$ and $p >71$ a prime and he proved that this equation has no non-trivial primitive integer solutions such that $10|xy$. 
In the same paper, he also proved the following result. 

\begin{theorem}
\begin{itemize}
\item[(a)] Suppose that $p>71$ is a prime. Then, for any $\alpha \in \{1,2,3,4\}$ and $\beta \in \{3, 4\}$, there are no non-trivial primitive integer solutions to the equation
$$
7x^p + 2^i 5^j y^p = 3z^5.
$$
\item[(b)] Suppose that $p>41$ is a prime. Then, for any $\alpha \in \{1,2,3,4\}$ and $\beta \in \{2, 3, 4\}$, there are no non-trivial primitive integer solutions to the equation
$$
x^5 + 7y^5 = 2^{\alpha} 5^{\beta} z^p.
$$
\end{itemize}
\end{theorem}

\section{Equations of signature $(p,q,r)$}
\label{sec:pqr}

This section discusses Eq.~\eqref{eq:genfermat} with exponents $p\neq q\neq r$. We will sometimes consider equations with some of the exponents being composite.

\subsection{The case $a=c=1$}

In this case, Eq.~\eqref{eq:genfermat} reduces to
\begin{equation}
\label{eq:xppbyqmzr}
 x^p+by^q=z^r.
\end{equation}

In 2014, Zhang and Luo~\cite{MR3222262} studied Eq.~\eqref{eq:xppbyqmzr} for $p=2$,
$r=3$ and $q$ even, that is,
\begin{equation}
\label{eq:x2pby2qmz3}
 x^2+by^{2q}=z^3
\end{equation}
and they obtained the following results. 

\begin{theorem}
Let $q \geq 7$ be prime and $\alpha \geq 6$ and $\beta \geq 0$ be integers. If $b \in \{ 2^{\alpha},4 \cdot 3^{2\beta}\}$, then \eqref{eq:x2pby2qmz3} has no non-trivial solutions in pairwise coprime integers $x,y,z$, except $(x,y,z,b)=(\pm 11, \pm 1,5,4)$.
\end{theorem}

\begin{theorem}
\label{th2:MR3222262}
 Let $q \geq 11$ be prime and let $\alpha \geq 6$, $\beta \geq 1$ and $\gamma \geq 0$ be integers that do not satisfy any of the following conditions. 
\begin{itemize}
  \item[(i)] {$2 \vert  \alpha$}, $2 \vert  \beta$, $2 \nmid \gamma$,
  \item[(ii)] {$2 \vert  \alpha$}, $2 \nmid \beta$, $2 \vert  \gamma$,
  \item[(iii)] $2 \nmid \alpha \beta \gamma$.
\end{itemize}
 If $b= 2^\alpha 3^\beta 31^\gamma$ and $q \nmid b$, then Eq.~\eqref{eq:x2pby2qmz3} has no non-trivial solutions in pairwise coprime integers $x,y,z$. Also, if $\gamma =0$, then the same result holds for $q=7$. 
\end{theorem}

Noubissie~\cite{noubissie2020generalized} studied Eq.~\eqref{eq:xppbyqmzr} with $p$ even,
$q=2$ and $b=r=3$, that is,
\begin{equation}
\label{eq:x2np3y2pmz3}
 x^{2p}+3y^{2}=z^3,
\end{equation}
and they obtained the following result.
\begin{theorem}
Eq.~\eqref{eq:x2np3y2pmz3} with integer $p \geq 8$ has no non-trivial primitive integer solutions.
\end{theorem}

Noubissie~\cite{noubissie2020generalized} also proved that for prime $q$ satisfying
$q >1964$, the equation
\begin{equation*}
x^{2}+3y^{2q}=z^3
\end{equation*}
has no non-trivial primitive integer solutions with $q\vert y$, while equation
\begin{equation*}
x^2+2y^3=z^{3r}
\end{equation*}
with prime $r \geq 3$ has no non-trivial primitive solutions with $x$ even.

Eq.~\eqref{eq:xppbyqmzr} with $b$ a positive odd integer, $p$ even, and
$5 \nmid h(-b)$ where $h(-b)$ denotes the class number of $\mathbb{Q}(\sqrt{-b})$ was considered by Zhang
and Bai~\cite{zhang2013} who obtained the following partial result.
\begin{theorem}
Let $p \geq 7$ be a prime, $b$ be a positive odd integer such that $5 \nmid  h(-b)$. Then the equation
\begin{equation*}
 x^{2p}+by^2=z^5
\end{equation*}
 has no solutions in pairwise coprime integers with $xyz \neq 0$ and $x$ even.
\end{theorem}

In the same paper, Zhang and Bai~\cite{zhang2013} considered Eq.~\eqref{eq:xppbyqmzr} with
$p=2$, $q$ even and $r=5$ and proved the following theorem.

\begin{theorem}
Let $\beta \geq 1$, $\gamma \geq 0$ and $p\geq 11$ prime. Let $b$ satisfy one of the following conditions.
\begin{itemize}
  \item[(i)] {$b=2^\alpha$}, $\alpha \geq 2$,
  \item[(ii)] {$b=2^\alpha 5^\beta 11^\gamma$}, $\alpha \geq 3$, $p \nmid b$,
  \item[(iii)] {$b=2^\alpha 5^\beta 19^\gamma$}, $\alpha \geq 2$, $p \nmid b$,
  \item[(iv)] {$b=2^\alpha 5^\beta 23^\gamma$}, $\alpha \geq 2$, $2 \vert  \alpha \beta \gamma$, $p \nmid b$, $\gamma \geq 1$ and $p \geq 13$.
\end{itemize}
 Then the equation
\begin{equation*}
 x^2+by^{2p}=z^5
\end{equation*}
 has no integer solutions in pairwise coprime integers $x,y,z$ with $xyz\neq 0$. The same conclusion holds in case (i) for $p=7$.  
\end{theorem}

The authors of~\cite{MR3222262} also considered Eq.~\eqref{eq:xppbyqmzr} with $q=2p$ and
$r=2$, that is,
\begin{equation}
\label{eq:xppby2pmz2}
 x^p+by^{2p}=z^2,
\end{equation}
and they obtained the following result.

\begin{theorem}
Let $p \geq 7$ be prime, and let $b=4 \cdot 3^{2\beta +1}$ with $\beta \geq 0$. Then Eq.~\eqref{eq:xppby2pmz2} has no non-trivial solutions in pairwise coprime integers $x,y,z$.
\end{theorem}

Noubissie~\cite{noubissie2020generalized} studied Eq.~\eqref{eq:xppby2pmz2} with $b=3^{2p-6}$, that
is,
\begin{equation*}
x^p+3^{2p-6}y^{2p}=z^2,
\end{equation*}
and proved that for prime $p \equiv 13$ modulo $24$, the equation has no non-trivial primitive integer solutions with $z \equiv 2(\text{mod}\, 4)$.

\subsection{The general $a,b,c$}

In this section we discuss cases of Eq.~\eqref{eq:genfermat} with non-equal exponents $(p,q,r)$ such that at least two of the coefficients $a,b,c$ are greater than $1$. 

In 2022, Cha\l upka et al.~\cite{CHALUPKA2022153} studied equations
\begin{equation}
\label{eq:ax2py2pm4z3}
 x^{2p}+ay^2 = 4z^3
\end{equation}

\begin{equation}
\label{eq:x2pay2pm4z3}
 x^2 +ay^{2p} = 4z^3
\end{equation}
where $p \geq 2$ is an integer.

In the case $a=7$ in \eqref{eq:ax2py2pm4z3}, they obtained the following results.

\begin{theorem}[Theorem 1 and 12 of~\cite{CHALUPKA2022153}]
 Let $p$ be a prime satisfying $5 \leq p <10^9$ with $p \neq 7,13$, then the equation
\begin{equation}
\label{eq:7x2py2pm4z3}
  x^{2p} +7y^2= 4z^3
\end{equation} 
 has no non-trivial primitive integer solutions.
\end{theorem}

\begin{theorem}
Let $p$ be a prime satisfying
\begin{equation*}
 p \equiv 3 \text{ or } 55 (\text{mod }106) \quad\text{ or } \quad p \equiv 47, 65, 113, 139, 143 \text{ or } 167  (\text{mod }168).
\end{equation*}
 Then Eq.~\eqref{eq:7x2py2pm4z3} has no non-trivial primitive integer solutions.
\end{theorem}

\begin{theorem}
If $a \in \{11, 19, 43, 67, 163\}$ and $p \geq 2$ is an integer, then Eqs.~\eqref{eq:ax2py2pm4z3} and \eqref{eq:x2pay2pm4z3} have no non-trivial primitive integer solutions.
\end{theorem}

Noubissie~\cite{noubissie2020generalized} also considered Eqs.~\eqref{eq:ax2py2pm4z3} and \eqref{eq:x2pay2pm4z3}
with $a=3$, and proved that for prime $p \geq 3$ these equations have no non-trivial
primitive integer solutions with $6 \vert y$. 

Cha\l upka et al.~\cite{CHALUPKA2022153} also studied some Eqs.~\eqref{eq:genfermat} with
$p=2q$. Specifically, they proved the following results.

\begin{theorem}
Let $p$ be a prime with $p \equiv 5$ (mod $6$). Then the equation
\begin{equation}
\label{eq:x2pm4ypm21z2}
  x^{2p} + 4y^{p} = 21z^2
\end{equation}
 has no non-trivial solutions.
\end{theorem}

They also prove that for all primes $p \geq 11$ with $p \neq 19$, Eq.~\eqref{eq:x2pm4ypm21z2} has no solutions in coprime odd integers $(x,y,z)$.

\begin{theorem}
Let $11 \leq p < 10^9$ and $p \neq 13, 17$ be a prime. Then the equation
\begin{equation}
\label{eq:32pm3x2pm4ypm7z2}
  3^{2p-3}x^{2p}+4y^p = 7z^2
\end{equation}
 has no primitive solutions $(x,y,z)$ in odd integers.
\end{theorem}

\begin{theorem}
Let $p$ be a prime satisfying
\begin{equation*}
 p \equiv 3 \text{ or } 55  (\text{mod }106) \quad\text{ or } \quad p \equiv 47, 65, 113, 139, 143 \text{ or } 167 (\text{mod }168).
\end{equation*}
 Then Eq.~\eqref{eq:32pm3x2pm4ypm7z2} has no primitive solutions in odd integers.
\end{theorem}

\begin{theorem}
The equation
\begin{equation*}
 x^{38} + 4y^{19} = 21z^2
\end{equation*}
 has no solution in coprime odd integers $(x,y,z)$.
\end{theorem}

In 2024, Cazorla Garc\'\i a~\cite{MR4793291} studied equation
\begin{equation}
\label{eq:garcia24}
 a x^2 + s^k y^{2n} = z^n
\end{equation}
where $s$ is prime, and provided an algorithmically testable set of conditions on $(a,s)$ which, if satisfied, imply the existence of a constant $B=B(a,s)$ such that if $n>B$ then for any $k\geq 0$ Eq.~\eqref{eq:garcia24} has no solutions in coprime integers $(x,y,z)$ such that $\gcd(ax, sy, z) = 1$ and $z$ is even. 

In a follow-on paper to~\cite{CHALUPKA2022153}, the authors~\cite{CHALUPKA2024} prove that for any
prime $p=5$ or $11$ modulo $12$, the equation
\begin{equation*}
3^{2p-3} 7 x^{2p}+4y^p=z^2
\end{equation*}
has no solutions in coprime odd integers.

In the same paper~\cite{CHALUPKA2024}, the authors prove the following result. 
\begin{theorem}
Let $n \geq 2$ be an integer. The equation
\begin{equation*}
 x^2+7y^{2n}=4z^{12}
\end{equation*}
 has no non-trivial primitive integer solutions $x,y,z$. 
\end{theorem}

In Chapter $5$ of their thesis, Putz~\cite{ThesisPutz} solved a number of equations of
signature $(5,3,11)$, that is, equation
\begin{equation}
\label{eq:pqr5311}
 ax^5+by^3=cz^{11}.
\end{equation}
Specifically, they proved the following theorem.
\begin{theorem}
Let
\begin{equation}
\label{eq:putz}
\begin{aligned}
   \mathcal{C}\coloneq  & \{(11^8,1,3^6 \cdot 5^{10}),(3^6 \cdot 11^8,1,5^{10}),(3^9 \cdot 11^8,1,5^{10}), \\ & (1,11^4,3^6 \cdot 5^{10}),  (3^4 \cdot 11^8,1,5^{10}),(11^8,3^3,5^{10}), \\ & (3^7 \cdot 11^8,1,5^{10}),   (3^8\cdot 11^8,1,5^{10}),(3^{10}\cdot 11^8,1,5^{10})\}.
\end{aligned}
\end{equation}
 If $(a,b,c) \in \mathcal{C}$ then Eq.~\eqref{eq:pqr5311} has no non-trivial primitive integer solutions.
\end{theorem}

\section{Summarizing table}
\label{sec:summary}

Table \ref{tbl10} below summarizes the equations solved in the literature considered in this survey. The table lists Eqs.~\eqref{eq:genfermat} with $|abc|>1$, see Table \ref{tb:pqrsigs} for the case $|abc|=1$. Table \ref{tbl10} specifies the equation, any restrictions on the parameters, the progress towards the solution of the equation, and the reference to the solution. The progress we refer to can be:

\begin{itemize}
 \item ``Solved'' implies that all primitive integer solutions have been described. This includes where the equation has no non-trivial primitive integer solutions. 
 \item ``Partial'' implies that there is only partial progress in the solution of the equation. This means that only the primitive solutions satisfying some additional conditions have been described.
 \item ``Non-trivial'' implies that it has been determined whether the equation has any non-trivial integer solution (that is, one with all variables non-zero).
 \item ``Non-zero'' implies that it has been determined whether the equation has any integer solution with $(x,y,z)\neq (0,0,0)$.
\end{itemize}

\captionof{table}{\label{tbl10} Generalised Fermat Eqs.~\eqref{eq:genfermat} with $1/p+1/q+1/r<1$ and $|abc|>1$ solved in the literature. Assume $m,n$ are integers and $p,r,s$ are primes. Unless otherwise stated, assume $\alpha, \beta, \gamma \geq 0$.}	
\begin{longtable}[c]{ |c|c|c|c|c|c|c|c|c| } 
 \hline
 Equation & Parameters  &  Progress & References \\
\hline \hline \endfirsthead 
 \multicolumn{4} {|c|}{\textbf{Special Signatures}}\\\hline \hline   \multicolumn{4} {|c|}{Signature $(4,4,4)$}\\\hline   \hline $x^4+y^4=cz^4$ & $1 \leq c \leq 100$, $c\neq 82$ & Solved & \cite{flynn2001covering,puadurariu2023rational} \\\hline
 $x^4+y^4=cz^4$ & $c \leq 10000$ & Non-zero & \cite{MR709852,MR4458887} \\ &&& \cite[Sec 6.6]{MR2312337} \\\hline
 $x^4+y^4=cz^4$ & $c >2$ fourth-power free,  & Solved & \cite{Grigorov1998Heights} \\ & rank of $v^2=u^3-cu$ is 1 & & \\\hline
 $x^4+2y^4=3z^4$ & & Solved &  Mathoverflow\footnotemark \\\hline 
 $x^4+3y^4=4z^4$ & &Solved &  \cite[pp. 630]{dickson1920history} \\\hline
 $x^4+8y^4=9z^4$ & &Solved&  \cite[pp. 630]{dickson1920history} \\\hline
 $ax^4+by^4=cz^4$ & $| a| +| b| +| c|  \leq 62$ & Non-zero & \cite[Sec 6.2.2]{mainbook} \\  &&& \cite[Sec 6.4-5]{wilcox2024systematic} \\\hline 
 $ax^4+by^4=bz^4$ & $| a| +2| b|  \leq 38$ & Non-trivial & \cite[Sec 6.15]{wilcox2024systematic} \\ \hline
 $x^4+b y^4=z^4$ & $1\leq b \leq 218$ &Non-trivial & \cite[Sec 6.15]{wilcox2024systematic} \\\hline 
 $x^4+34y^4=z^4$ &&Solved& \cite{soderlund2020note} \\\hline
 $x^4+2py^4=z^4$ & Any of the conditions  & Solved &  \cite{taclay2023} \\ & in Theorem \ref{th:x4p2py4mz4} && \\\hline 
 $a^3x^4+b^3y^4=c^3z^4$ & Conditions in  & Solved & \cite{MR249355,Mordell_1970} \\  & Theorem \ref{th:MR249355} && \\\hline \hline
 \multicolumn{4} {|c|}{Signature $(2,4,r)$}\\\hline  \hline  $2x^2+y^4=z^n$ & $n \geq 5$ & Solved & \cite{bennett2010diophantine} \\\hline
 $ax^2 +y^4 = z^p$ & $(a=3,p > 17)$ or & Solved & \cite{dieulefait2008solvingfermattypeequations,MR4473105} \\ &  $(a=5,p > 499)$ or && \\ & $(a=6,p > 563)$ or  &&  \\ & $(a=7,p > 349)$ && \\\hline  
 $sx^p+y^4=z^2$ & $s >3, s \equiv 3$ (mod 8), &Solved & \cite{MR2737959} \\ &  $s \neq 2t^2+1$, && \\ &  $p > (8\sqrt{s+1}+1)^{16(s-1)}$ && \\\hline
 $sx^p+y^4=z^2$ & $s \equiv 5$ (mod 8), $s \neq t^2+4$, &Solved & \cite{MR2737959} \\ &  $p > (8\sqrt{s+1}+1)^{16(s-1)}$ && \\\hline
 $2^{\alpha}x^p+y^4=z^4$ & $\alpha \geq 0$, $p \geq 5$ & Solved& \cite{MR2737959} \\\hline
 $2^{\alpha}s^{\beta}x^p+y^4=z^4$ & $\alpha \geq 0$, $\beta \geq 0$, $s\geq 3$,& Solved& \cite{MR2737959,MR4205757} \\  &  $s \neq 2^{2^k}+ 1$, && \\ &  $p > (\sqrt{8s+8}+1)^{2s-2}$ && \\\hline
 $sx^p+y^4=z^4$ & $s \in \{5,17\}$, $p>5$ & Solved & \cite{MR4205757} \\\hline
 $2^{\alpha}s^{\beta}x^p+y^4=z^4$ & $2 \leq s <50$, $s \neq 5$, $s \neq 17$,   & Solved & \cite{MR4205757} \\  & $p>5$, $\alpha, \beta \geq 0$ && \\\hline 
 $x^p+py^2=z^4$ & $p \equiv 1$ (mod 4), $p \nmid B_{(p-1)/2}$ & Partial & \cite{MR1341665} \\\hline  
 $ax^{2m}+y^2=z^n$ & Conditions in & Solved & \cite{MR1604052} \\ &  Theorem \ref{th:ax2mpy2mzn} && \\\hline
 $x^p+by^2=z^4$ & $p >19$, $p \equiv 1,3$ (mod 8),  & Solved & \cite{MR4609012} \\ & $(b=6,p \neq 97)$ or  && \\ & $(b=10, p \neq 139)$ or   && \\ & $(b=11,p \neq 73)$ or   && \\ & $(b=19,p\neq 43,113)$ or && \\ & $(b=129,p \neq 43)$ && \\\hline
 $x^p+129y^2=z^4$ & $p >900$  & Solved& \cite{MR4609012} \\ \hline \hline
 \multicolumn{4} {|c|}{Signature $(2,6,r)$}\\\hline \hline  $x^p +2y^6= z^2$ & $p \neq 7$, $p \equiv 1,7$ (mod 24) & Solved & \cite{Chen_2012} \\\hline  
 $x^2 +2y^6= z^p$ & $p >257$ & Solved & \cite{MR4473105} \\\hline 
 $x^2 +3y^6= z^n$ & $n \geq 3$ & Solved & \cite{koutsianas2019generalizedfermatequationa23b6cn} \\\hline 
 $x^2 +6y^6= z^p$ & $p >563$ & Solved & \cite{MR4473105} \\\hline 
 
 $x^2+by^6=z^p$ & $p \geq r_0$, $(b,r_0)$ given by & Solved &  \cite{MR4583916} \\ &  Table \ref{tb:x2pby6mzr} &&  \\\hline
 $x^2+b^2y^{2m}=z^{4n}$ & $b>0$ odd with prime &Partial & \cite{MR1188732} \\  &   factor $l = \pm 3$ (mod 8)  && \\\hline
 \hline \multicolumn{4} {|c|}{Signature $(3,3,4)$}\\\hline  \hline $2x^3+2y^3=z^4$ & & Solved & \cite[App. A]{wilcox2024systematic} \\\hline \hline
 \multicolumn{4} {|c|}{Signature $(4,4,3)$}\\\hline \hline  $x^4+2y^3=z^4$ & & Solved & \cite[App. A]{wilcox2024systematic} \\\hline
 $x^4+3y^3=z^4$ & & Solved & \cite[App. A]{wilcox2024systematic} \\\hline
 $x^4+2y^4=z^3$ & & Solved & \cite{soderlund2017primitive} \\\hline
 $x^4+3y^4=z^3$ & & Solved & \cite{MR1288426} \\\hline
 $x^4+s^2y^4=z^3$ & & Solved & \cite{soderlund2020diophantine} \\\hline
 $x^4+by^4=z^3$ & $b \in \{19,43,67\}$ & Solved & \cite{soderlund2019some} \\\hline \hline
 \multicolumn{4} {|c|}{\textbf{Signature $(p,p,p)$}}\\\hline  \hline \multicolumn{4} {|c|}{Explicit small values of $p$}\\\hline \hline   $x^5+y^5=cz^5$ & $1 \leq c \leq 100$, $c \neq 88$ & Non-trivial &\cite[Sec 6.3.4]{mainbook} \\\hline
 $x^5+y^5=cz^5$ & $c$ only divisible by primes & Solved & \cite{kraus2004} \\ & $c_i \not \equiv 1$ (mod 5) && \\\hline
 $ax^5+by^5=cz^5$ & $a+b+c< 19$ &Non-zero& \cite[Sec 6.2.3]{mainbook} \\  &&& \cite[Sec 6.7]{wilcox2024systematic} \\\hline 
 $x^6+y^6=cz^6$ & $c \leq 164634913$ & Non-trivial &  \cite{MR4552507} \\\hline
 $ax^6+by^6=bz^6$ & $| a| +2| b|  \leq 10$ &Non-trivial&  \cite[Sec 6.16]{wilcox2024systematic} \\\hline
 $x^{14}+2^{\alpha} 7^{\beta} y^{14}=z^{14}$ & $\alpha \geq 0,\beta \geq 1$ & Solved & \cite[pp. 736]{dickson1920history}  \\\hline
 \hline
 \multicolumn{4} {|c|}{ The case $a=b=1$}\\\hline \hline   $x^p+y^p=15z^p$ & $p \geq 5$ and $2p+1$ or & Solved & \cite{kraus1996equations} \\ &  $4p+1$ is prime && \\\hline
 $x^p+y^p=15z^p$ & $ p \geq 5$ & Partial &  \cite{kraus1996equations} \\\hline 
 $x^p+y^p=pz^p$ & $p$ regular & Solved & \cite[pp. 759]{dickson1920history} \\ &  (see Theorem \ref{th:regular}) && \\\hline
 $x^p+y^p=pz^p$ & $p\geq 3$ & Partial & \cite{MR1517149}  \\\hline
 $x^n+y^n=2^{\alpha}z^n$ & $\alpha \geq 0$, $n \geq 3$ & Solved & \cite{Darmon1997,ribet1997equation} \\\hline 
 
 $x^p+y^p=31^\alpha z^p$ & $7  \leq p \leq 10^6$,  $\alpha \geq 0$   & Solved & \cite[Sec 15.7]{MR2312338} \\ \hline
 
 $x^p+y^p=s^\alpha z^p$ & $s\neq 2^n\pm 1$ odd, $\alpha \geq 0$, & Solved & \cite{MR1611640} \\ & $p > \left(1+\sqrt{(s+1)/2}\right)^{(s+11)/6}$ && \\\hline
 $x^p+y^p=s^{\alpha}z^p$ & $3 \leq s \leq 100$, $p \geq 5$,  & Solved & \cite{ribet1997equation,MR1611640}, \\& $\alpha \geq 0$, $s \neq 31$ && \cite[Th.15.5.3]{MR2312338} \\\hline
 $x^p+y^p=psz^p$ & $p \geq 5$ regular,  & Solved & \cite{sitaraman2000fermat} \\ &  $s$ only divisible by primes && \\ & $kp-1$ with $\gcd(k,p)=1$ && \\\hline \hline
 \multicolumn{4} {|c|}{ The case $a=1$}\\\hline  \hline $x^p+3y^p=5z^p$ & $p \geq 5$ and $2p+1$  & Solved & \cite{kraus1996equations} \\ & or $4p+1$ is prime && \\\hline
 $x^p+3y^p=5z^p$ & $p \geq 5$ & Partial & \cite{kraus1996equations} \\\hline
 $x^p+3y^p=5z^p$ & $5 \leq p < 10^7$ & Solved & \cite{Kraus2002} \\\hline
 $x^p+2^{\beta} y^p= s^{\alpha}z^p$ & $s > 3$, $s\neq 2^n \pm 1$,  & Solved & \cite{MR1611640} \\  & $p >(1+\sqrt{8(s+1)})^{2(s-1)}$,  && \\ & $\alpha,\beta \geq 0$, $p \nmid \alpha$  \\\hline
 $x^p+16y^p=s^{\alpha}z^p$ & $\alpha \geq 0$, $s\geq 3$, $s \neq 17$,  & Solved & \cite{MR1611640} \\   &   $p \neq s$, $p \geq 5$, && \\ & $p > \left(1+\sqrt{(s+1)/6}\right)^{(s+1)/6}$  && \\  \hline \hline
 \multicolumn{4} {|c|}{ General $a,b,c$}\\\hline  \hline $ax^n+by^n=16z^n$ & $(a,b) \in \{(25,23^4), (5^8,37), $ & Solved & \cite{Dahmen_2012} \\  & $(5^7,59^7),(7,47^7),(11,85^2)  \}, $  && \\ & $n \geq 5$ odd && \\\hline 
 $3x^p+4y^p=5z^p$ & $p \equiv 5$ (mod 8) or & Solved & \cite{Kraus2002,FREITAS2016751} \\  & $p \equiv 19$ (mod 24) &&  \\\hline
 $3x^p+8y^p=21z^p$ & $p \equiv 5$ (mod 8) or & Solved & \cite{FREITAS2016751} \\  & $p \equiv 23$ (mod 24) && \\\hline  
 $ax^p+by^p=16cz^p$ & Prime factors of $abc$  & Solved & \cite{MR4203704} \\  &  are 1 (mod 3), && \\ &  $p$ sufficiently large && \\\hline
 $ax^p+by^p=16cz^p$ & $n>0$ not a divisor of  & Solved & \cite{MR4203704} \\  &   14,16 or 18, prime factors  && \\ &  of $abc$ are $\pm 1$ (mod $n$), && \\ & $p$ sufficiently large && \\\hline
 $ax^p+by^p=2^rcz^p$ & $r \geq 0$, $r \neq 1$, & Solved & \cite{MR4203704} \\   & prime factors of $abc$  && \\ & are 1 (mod 12), && \\ &  $p$ sufficiently large &&  \\\hline
 $ax^p+by^p=cz^p$ & Conditions in  & Solved & \cite{MR4203704} \\  & Theorem \ref{th:qoddprime} && \\\hline  
 $ax^p+by^p+2^ncz^p=0$ & Any of the conditions & Solved & \cite{MR4203704} \\  &  in Theorem \ref{th:jacobi} && \\\hline
 $ax^p+2^rby^p+2^rcz^p=0$ & Any of the conditions & Solved & \cite{MR4203704} \\  &  in Theorem \ref{th:jacobi} &&  \\\hline
 $ax^n+by^n=cz^n$ &Conditions in & Solved & \cite{POWELL198434,MR779375} \\  &  Theorem \ref{th:powell}  && \\\hline \hline
 \multicolumn{4} {|c|}{\textbf{Signature $(p,p,r)$} }\\\hline  \hline \multicolumn{4} {|c|}{ The case $a=b=1$ and fixed $p$}\\\hline \hline   $x^3+y^3=2z^{2r}$ & $r \geq 2$ & Solved & \cite{zhang2014} \\\hline
 $x^3+y^3=s^\alpha z^r$ & $\alpha \geq 1$, $r \geq s^{2s}$, $s\geq 5$,  & Solved & \cite{MR3830208} \\ &  $s\not\in S_0$, see \eqref{eq:S0def} && \\\hline
   $x^3+y^3=s^{\alpha}z^r$ & $s=2^a 3^b-1$, $a \geq 5$,  & Solved & \cite{MR3830208} \\  &  $b \geq 1$, $\alpha \geq 1$, $r>s^{2s}$, && \\ & $(\alpha,r,a,b)$ do not && \\ &  satisfy \eqref{cd:th2:bennett2017} && \\\hline
 $x^3+y^3=s z^r$ & $s\not\in T$, see Theorem \ref{th3:bennett2017}, & Solved & \cite{MR3830208} \\  & positive proportion of $r$ && \\\hline 
 $x^3+y^3=sz^r$ & $(s,r)$ satisfy Table \ref{tb:MR3830208},  & Solved & \cite{MR3830208}   \\ & $r \geq s^{2s}$ && \\ \hline
 $x^3+y^3=5z^r$ & $r\equiv 13,19$ or 23 (mod 24) & Solved & \cite{MR3830208}   \\ \hline
 $x^4+y^4=sz^r$ & $s \in \{73,89,113 \}$, $r >13$ & Solved & \cite{MR2139003} \\\hline  
 $x^5+y^5=s^{\alpha}z^r$ & $r \geq s^{13s}$, $\alpha \geq 1$,  & Solved & \cite{bruni2015twisted} \\ & $s \not \in {\cal{W}}_5$ prime, && \\ & see Theorem \ref{th:bruni2015}  && \\\hline  
 $x^5+y^5=2z^r$ & & Partial & \cite{signature2019multi} \\\hline
 
 $x^p+y^p=3z^n$ & $p \in \{5,7, 13 \}$, $n \geq 2$ & Solved & \cite{signature2019multi,10.5565/PUBLMAT6722309,billerey2024darmonsprogramgeneralizedfermat} \\\hline 

 $x^5+y^5=2c z^r$ & $c$ only divisible by primes & Solved & \cite{44370229} \\ & $c_i \not \equiv 1$ (mod 5),  $r > 13$,  && \\ & $r \equiv 1$ (mod 4) or && \\ & $r \equiv \pm 1$ (mod 5) && \\\hline 
 $x^5+y^5=3c z^r$ & $c$ only divisible by primes & Solved & \cite{44370229} \\ & $c_i \not \equiv 1$ (mod 5),  $r > 73$,  && \\ & $r \equiv 1$ (mod 4) or && \\ & $r \equiv \pm 1$ (mod 5) && \\\hline
 $x^5+y^5=7z^r$ & $r \geq 13$ & Solved & \cite{billerey2008solvingfermattypeequationsx5y5dzp} \\\hline
 $x^5+y^5=13z^r$ & $r \geq 19$ & Solved & \cite{billerey2008solvingfermattypeequationsx5y5dzp} \\\hline 
 $x^5+y^5=2^{\alpha} 3^{\beta} 5^{\gamma} z^r$ & $r \geq 7$, $\alpha \geq 2$, $\beta, \gamma \geq 0$ & Partial & \cite{Billerey_2007,billerey2008solvingfermattypeequationsx5y5dzp}  \\\hline
 $x^5+y^5=2^{\alpha} 3^{\beta} 5^{\gamma} z^r$ & $r \geq 13$, $\alpha \geq 2$, $\beta, \gamma \geq 0$ & Solved & \cite{billerey2008solvingfermattypeequationsx5y5dzp} \\\hline
 $x^5+y^5=2^{\alpha}11z^r$ & $r > 19$, $\alpha \geq 2$ & Solved &\cite{noubissie2020generalized}  \\\hline
 $x^5+y^5=2^{\alpha}11z^r$ & $r > 19$, $\alpha=0,1$ & Partial &\cite{noubissie2020generalized}  \\\hline
 $x^6+y^6=2z^r$ & $r \geq 5$ & Solved & \cite{kraus2002question} \\\hline  
 $x^p-y^p=cz^r$ & $(p,r) \in \{(6,2),(8,2),(8,3), $  & Solved & \cite{kraus2002question} \\ & $(9,3),(10,2),(12,2)\}, $ && \\ & $c \in {\cal{D}}_p$, see Theorem \ref{th:kraus2002} && \\\hline
 $x^p-y^p=cz^2$ & $p\in\{5,7\}$, $p \vert c$,  & Solved & \cite{kraus2002question} \\ & $c \in {\cal{D}}_p$, see Theorem \ref{th:kraus2002} && \\\hline 
 $x^{12}-y^{12}=cz^n$ & $n \geq 3$, $c \geq 1$,   & Solved & \cite{kraus2002question} \\  & $l$ prime divisors of $c$, && \\ & $\nu_l(c)<12$ && \\ & and $l \not \equiv 1$ (mod 4) && \\\hline
 $x^9+y^9=2^nz^2$ & $n \geq 0$ & Solved & \cite{MR995897} \\\hline
 $x^{13}+y^{13}=s \gamma z^r$ & $s \in \{3,5,7,11 \}$,  & Partial & \cite{FreitasDieulefait2013} \\  & $\gamma$ divisible only by primes && \\ &  $\gamma_i \not \equiv 1$ (mod 13), && \\ &  $r \not \in R$, see \eqref{def:R} && \\\hline
 $x^{14}+y^{14}=2^{\alpha} 3^{\beta} 5^{\gamma} c z^r$ & $r \geq 17$, $\alpha \geq 2$ or  & Solved & \cite{freitas2015recipes} \\  & $\beta> 0$ or $\gamma > 0$ && \\ & $c$ only divisible by primes && \\ & $c_i \not \equiv 0,1$ (mod 7)  && \\\hline
 
 $x^{26}+y^{26}=10 \gamma z^r$ &$\gamma$ only divisible by primes & Solved & \cite{freitas2014criteriairreducibilitymodp} \\  &  $\gamma_i \not \equiv 1$ (mod 13), && \\ &  $r \not \in R$, see \eqref{def:R}  && \\\hline \hline
   \multicolumn{4} {|c|}{ The case $a=b=1$ and fixed $r$}\\\hline \hline  $x^5+y^5=2^{\alpha}5^{\beta}s^{\gamma}z^2$ & $s=\pm 3$ (mod 10),  & Solved & \cite{soderlund2011diophantine,soderlund2013note} \\  &  $\alpha \geq 0$, $\beta \geq 0$, $\gamma \geq 0$ && \\\hline
 $x^n+y^n=cz^2$ & $c \in \{ 2,3,5,6,10, $ & Solved & \cite{Bennett_Skinner_2004} \\ & $ 11,13,17 \}$, $n \geq 5$ && \\ \hline
 $x^p+y^p=cz^2$ & $c \geq 5$, $c \neq 7$, $p >c^{12c^2}$ & Partial & \cite{bennett2006diophantine} \\\hline 
 $x^p+y^p=2sz^2$ & $s \geq 5$, $p >s^{132s^2}$ & Solved & \cite{bennett2006diophantine} \\\hline
 $x^p+y^p=2^{\alpha}sz^2$ & $p>s^{27s^2}$,  $s \neq 7$, $\alpha \geq 1$ & Solved & \cite{mulholland2006elliptic} \\ \hline
 $x^p+y^p=2^{\alpha}cz^2$ & $c$ odd square-free, $3 \nmid c$,   & Solved & \cite{MR2737959} \\  &  $7 \nmid c$, $p>c^{132c^2}$, $\alpha \geq 1$ && \\\hline
 $x^p+y^p=cz^2$ & $c \geq 3$, $c$ square-free,  & Solved & \cite{MR2310336} \\ & $c$ only divisible by primes && \\ & $c_i \not \equiv 1$ (mod $p$),  && \\ & $p \in \{11,13,17\}$ && \\\hline
 $x^n+y^n=cz^3$ & $c \in \{2,3,4,5 \}$, $n\geq 4$ & Solved &\cite{bennett2004ternary}  \\\hline
 $x^p+y^p=s^{\alpha} z^3$ & $\alpha \geq 0$, $p >s^{4s^2}$ & Solved & \cite{bennett2004ternary} \\\hline
 $x^p+y^p=cz^3$ & $c \neq 0$ cube-free, $3 \vert c$,   & Solved & \cite{KarolinaKrawciow2011} \\ & $C=\mathrm{rad}(c)$, $p > C^{10C^2}$ &&  \\\hline \hline
 \multicolumn{4} {|c|}{ The case $a=c=1$ and fixed $r$}\\\hline \hline   $x^5+2y^5=z^2$ & & Solved & \cite{ivorra2003equations} \\\hline
 $x^p+2^{\alpha}y^p=z^2$ & $\alpha \geq 2$, $p \geq 5$ & Solved & \cite{ivorra2003equations} \\ &&& \cite[Sec 15.3]{MR2312338} \\\hline  
 $x^p+s^{\alpha}y^p=z^3$ & $\alpha \geq 0$, $p>s^{2s}$, & Solved & \cite{bennett2004ternary} \\ &  $s \neq s_1^3\pm 3^{s_2}$ for $s_2 \neq 1$  && \\\hline
 $x^p+s^{\alpha}y^p=z^3$ & $s \in \{ 5,11,13,23,29,31, $ & Solved & \cite{bennett2004ternary} \\ & $41,43,47,53,59,61,67, $ &&\\ & $71,79,83,97 \}$ && \\ & $\alpha \geq 1$, $p \geq 11$,  $p \nmid s^2-1$, && \\ & $(s, p) \not\in \{(13, 19),(29, 11), $ && \\ & $(43, 13),(47, 13),(59, 11), $ && \\ & $(61, 61),(67, 73),(79, 97), $  && \\ & $(97, 13),(97, 79)\}$ &&  \\\hline
 $x^4+by^4=z^5$ & $b \in \{3,11,19,43,67,163\}$ & Solved & \cite[Th.2.8]{soderlund2019some} \\\hline
 \hline
 \multicolumn{4} {|c|}{ The case $a=1$ and fixed $p$}\\\hline \hline
 $x^5 + 7 y^5 = 2^{\alpha} 5^{\beta} z^r$ & $r >41, \alpha \in \{1,2,3,4\}, \beta \in \{2,3,4\}$ & Solved & \cite{azon2025} \\\hline
 $x^5+7y^5=z^r$ & $r>41$ & Partial & \cite{azon2025} \\\hline \hline
 \multicolumn{4} {|c|}{ The case $a=1$ and fixed $r$}\\\hline  \hline $x^p+2^{\beta}y^p=2z^2$ & $p \geq 7$, $1 \leq \beta \leq p-1$ & Solved & \cite{ivorra2003equations} \\  &&& \\\hline
 $x^p+2^{\alpha}y^p=cz^2$ & $c \in \{3,5,7,11,13,15,17 \}$,  & Solved & \cite{Bennett_Skinner_2004} \\  & $\alpha \geq \alpha_0$, $p > c$, && \\ & $(c,\alpha,p) \neq (11,3,13)$, && \\ & where $(c,\alpha_0)$ && \\ & is given by \eqref{eq:BennettSkinner2004} && \\\hline 
 
 $x^p+s^{\alpha}y^p=2z^2$ & $p \geq 11$, $s \in \{5,11,13 \}$, & Solved & \cite{Bennett_Skinner_2004} \\  &  $s \neq p$, $\alpha \geq 0$ && \\\hline
 $x^p+2^\alpha y^p=sz^2$ & $s \neq (2^k \pm 1)/d^2$, $k,d \in \mathbb{Z}$,  & Solved & \cite{zhang2012} \\  & $\alpha \geq 0$, $\alpha \neq 1$, $p > s^{8s^2}$ && \\\hline
 $x^p+2^\alpha y^p=sz^2$ & $\alpha \geq 2$, $s >5$,  & Solved & \cite{zhang2012} \\  &  $s =\pm 3$ (mod 8),  && \\ & $p >s^{8s^2}$ && \\\hline
 $x^p + s^{\alpha}y^p = 3^{\beta}z^3$ & $s \not\in \{ 5, 3t^3 \pm1, $ & Solved & \cite{bennett2004ternary} \\  & $9t^3 \pm 1 : t \in \mathbb{N}\}$ && \\  &  $\alpha, \beta \geq 1$, $3 \nmid \beta$, $p>s^{28s}$ && \\\hline 

 $x^p+s^{\alpha}y^p=3^{\beta}z^3$ & $p \geq 7$, $s \in \{7,11,13\}$, & Solved & \cite{bennett2004ternary}  \\  & $(s,p) \neq (7,13),(13,7)$,  && \\ & $3 \nmid \beta$, $\alpha,\beta \geq 1$ && \\\hline
 
 $x^p+3^{\alpha} y^p=c^{\beta} z^3$ & $c \in \{2,3,5,7,11,13, $  & Solved & \cite{bennett2004ternary} \\
  & $15,17,19\}$, && \\  & $(c,p) \not \in \{(7,11),(11,13)\}$, && \\ & $(\alpha,c) \not \in\{(1,t):t=2, \text{ or }$ && \\ & $ t \geq 11\}, $ && \\ &  $\alpha, \beta \geq 0$,  $p>\max\{c,4\}$ &&  \\\hline
 $x^p+s^{\alpha}y^p=Mz^3$ & Conditions of  & Solved & \cite{KarolinaKrawciow2011} \\  & Theorem \ref{th:KarolinaKrawciow2011} && \\\hline  
 \hline
 \multicolumn{4} {|c|}{ The case $c=1$ and fixed $p$ }\\\hline  \hline
 $ax^3 + by^3 = z^n$ & Conditions of  & Partial & \cite[Th.13.2]{bennett2013klein} \\ & Theorem \ref{th:bennett2013klein} && \\\hline 
 \hline
 \multicolumn{4} {|c|}{ The case $c=1$ and fixed $r$ }\\\hline   \hline
 $ax^p+by^p=z^2$& $p \geq 11$, $ab \in \{2^{\alpha}11^{\beta}, 2^{\alpha}13^{\beta}, $ & Solved & \cite{Bennett_Skinner_2004} \\ 
 & $2^{\alpha}19^{\beta},2^{\alpha}29^{\beta},2^{\alpha}43^{\beta}, $   && \\
 & $2^{\alpha}53^{\beta},2^{\alpha}59^{\beta},2^{\alpha}61^{\beta},2^{\alpha}67^{\beta}  \}$, && \\ & $a,b$ coprime, $p \nmid ab$, && \\ & $\alpha = 0$ or $\alpha \geq 3$, $\beta \geq 1$, && \\
 &  $ab,p,\alpha,\beta$ not in Table \ref{tb:ABvalsodd} &&  \\ \hline
 $ax^p+by^p=z^2$& $p \geq 11$, $ab \in \{ 2 \cdot 19^{\beta}, 4 \cdot 11^{\beta},  $ &Solved & \cite{Bennett_Skinner_2004} \\ & $4 \cdot 19^{\beta},4 \cdot 43^{\beta}, 4 \cdot 59^{\beta}, $ && \\ & $4\cdot 61^{\beta}, 4\cdot  67^{\beta} \}, $ && \\ & $a,b$ coprime, $p \nmid ab$,  $\beta \geq 1$,  && \\ 
 & $ab,p,\alpha,\beta$ not in Table \ref{tb:ABvalsodd} && \\ \hline
 
 $ax^p+by^p=z^2$& $a,b$ coprime, $ab=2^{\alpha}s^{\beta}t^{\gamma}$,   &  Solved  & \cite{Bennett_Skinner_2004} \\ &   $\beta,\gamma \geq 0$, $\alpha \geq 6$, && \\ & $(s,t)  \in \{   (3, 31) (\beta \geq 1), $ && \\ & $(5, 11) (\alpha \geq 7), (5, 19), $ && \\ & $(5, 23) (\beta \geq 1),  (7, 19) (\gamma \geq 1), $ && \\ & $  (11, 13),(11, 23) (\beta \geq 1), $ && \\ & $ (11, 29),(11, 31) (\beta \geq 1), $  && \\ & $ (13, 31) (\beta \geq 1),(19, 23) (\beta \geq 1)$, && \\ & $(19, 29),(29, 31) (\beta \geq 1) \}, $ && \\ & $p \geq 11$, $p \nmid ab$, && \\ & $ab,p,\alpha,\beta$ not in Table \ref{tb:ABvalseven} && 
 \\\hline
 $ax^p+by^p=z^2$ & $ab=2^\alpha 3^\beta$, $\gcd(a,b)=1$,  & Solved & \cite{MR3222262} \\ & $\alpha \geq 6$,  $\beta \geq 0$, $p \geq 7$, $p \nmid ab$  && \\\hline \hline
 \multicolumn{4} {|c|}{The general $a,b,c$ and fixed $r$ }\\\hline \hline   $5^{\alpha} x^p+64y^p=3z^2$ & $p \geq 7$, $\alpha \geq 1$ & Partial & \cite{armandphdthesis,Noubissie21} \\ \hline
 $ax^p+by^p=cz^2$ & $(a,b,c) \in\{ (2^n,l^m,1), $ & Solved & \cite{ivorra2006quelques} \\ & $(2^n l^m, 1,1),(1,l^m,2) \}$, && \\ & conditions below \eqref{abc_in_axppbypmcz2}  && \\\hline
 $4x^p+y^p=3z^2$ & $p \geq 7$ & Solved & \cite{ivorra2006quelques} \\\hline
 $64x^p+y^p=7z^2$ & $p \geq 11$ & Solved & \cite{ivorra2006quelques} \\  \hline
 $3^{31}x^{34}-4y^{34}=7z^2$ & & Partial & \cite{CHALUPKA2022153} \\  \hline
 $2^{\alpha}x^p+27y^p=s^{\beta} z^3$ & $s \in \{7,13\}$, $p >s $, $\alpha \geq 1$, & Solved & \cite{armandphdthesis,Noubissie21}\\  & $ \beta \geq 1$, $\beta \equiv 1$ (mod 3) && \\\hline
 $7x^p+y^p=3z^5$ & $p>71$ & Partial & \cite{azon2025} \\\hline
 $7x^p + 2^{\alpha} 5^{\beta} y^p = 3z^5$ & $p>71, \alpha \in \{1,2,3,4\}, \beta \in \{3,4\}$ & Solved & \cite{azon2025} \\\hline \hline
 \multicolumn{4} {|c|}{ \textbf{Signature} $(p,q,r)$ }\\\hline \hline   \multicolumn{4} {|c|}{The case $a=c=1$ }\\\hline  \hline $x^{2n}+3y^2=z^3$ & $n \geq 8$ & Solved & \cite{noubissie2020generalized} \\\hline
 $x^{2}+3y^{2p}=z^3$ & $p >1964$ & Partial & \cite{noubissie2020generalized} \\\hline   
 $x^2+2y^3=z^{3p}$ & $p \geq 3$ & Partial & \cite{noubissie2020generalized} \\\hline
 $x^2+by^{2p}=z^3$ & $p \geq 7$, $\alpha \geq 6$, $\beta \geq 0$,  & Solved & \cite{MR3222262} \\ & $b \in \{2^\alpha,4 \cdot 3^{2\beta}\}$ && \\\hline
 $x^2+2^\alpha 3^\beta 31^\gamma y^{2p}=z^3$ & Conditions of  & Solved & \cite{MR3222262} \\  &Theorem \ref{th2:MR3222262} && \\\hline
 $x^{2p}+by^2=z^5$ & $p \geq 7$, $b$ positive odd,  & Partial & \cite{zhang2013} \\ & $5 \nmid h(-b)$ &&  \\\hline
 $x^2+2^{\alpha}y^{2p}=z^5$ & $\alpha \geq 2$, $p \geq 7$ & Solved & \cite{zhang2013} \\\hline
 $x^2+2^{\alpha}5^{\beta}11^{\gamma}y^{2p}=z^5$ & $\alpha \geq 3$, $p >11$,  & Solved & \cite{zhang2013} \\  & $\beta \geq 1$, $\gamma \geq 0$ && \\\hline
 $x^2+2^{\alpha}5^{\beta}19^{\gamma}y^{2p}=z^5$ & $\alpha \geq 2$, $p \geq11$, $p \neq 19$,  & Solved & \cite{zhang2013} \\  & $\beta \geq 1$, $\gamma \geq 0$ && \\ \hline
 $x^2+2^{\alpha}5^{\beta}23^{\gamma}y^{2p}=z^5$ &  $\alpha \geq 2$, $\beta \geq 1$, $\gamma \geq 1$, $2\vert \alpha \beta \gamma$ & Solved & \cite{zhang2013} \\  & $p \geq 13$, $p \neq 23$ && \\\hline
 $x^p+4 \cdot 3^{2\beta+1}y^{2p}=z^2$ & $p \geq 7$, $\beta \geq 0$ & Solved & \cite{MR3222262} \\ \hline
 $x^p+3^{2p-6}y^{2p}=z^2$ & $p \geq 7$,  & Partial & \cite{noubissie2020generalized} \\  & $p \equiv 13$ (mod 24) && \\\hline
 \hline \multicolumn{4} {|c|}{The general $a,b,c$ }\\\hline  \hline $x^2+3y^{2p}=4z^3$ & $p \geq 3$ & Partial & \cite{noubissie2020generalized} \\\hline
 $x^{2p}+3y^{2}=4z^3$ & $p\geq 3$ & Partial & \cite{noubissie2020generalized} \\\hline
 $7x^2+y^{2p}=4z^3$ & $5 \leq p <10^9$, $p \neq 7,13$ & Solved & \cite{CHALUPKA2022153} \\\hline
 $7x^2+y^{2p}=4z^3$ & $p\equiv 3,55$ (mod 106) or & Solved & \cite{CHALUPKA2022153} \\ &  $p \equiv 47,65,113, $ && \\ & 139,143 or 167 (mod 168)&& \\\hline
 $ax^2+y^{2n}=4z^3$ & $a \in \{11,19,43,67,163\}$,  & Solved & \cite{CHALUPKA2022153} \\ & $n \geq 2$ && \\\hline
 $x^2+ay^{2n}=4z^3$ & $a \in \{11,19,43,67,163\}$, & Solved & \cite{CHALUPKA2022153} \\ &  $n \geq 2$ && \\\hline
 $x^{2p}+4y^p=21z^2$ & $p \equiv 5$ (mod 6) & Solved & \cite{CHALUPKA2022153} \\\hline
 $x^{2p}+4y^p=21z^2$ & $p \geq 11$ & Partial & \cite{CHALUPKA2022153} \\\hline
 $3^{2p-3}x^{2p}+4y^p=7z^2$ & $11 \leq p<10^9$, $p \neq 13,17$ & Partial & \cite{CHALUPKA2022153} \\ \hline
 $3^{2p-3}x^{2p}+4y^p=7z^2$ & $p\equiv 3,55$ (mod 106) or & Solved & \cite{CHALUPKA2022153} \\  &  $p \equiv 47,65,113, $ && \\ & 139,143 or 167 (mod 168)&& \\\hline
 $x^{38}+4y^{19}=21z^2$ & & Partial& \cite{CHALUPKA2022153} \\ \hline
   $ax^2+s^ky^{2n}=z^n$ & See conditions below \eqref{eq:garcia24} & Partial & \cite{MR4793291} \\ \hline
 $3^{2p-3}7x^{2p}+4y^p=z^2$ & & Partial & \cite{CHALUPKA2024} \\   \hline
 $x^2+7y^{2n}=4z^{12}$ & $n \geq 2$ & Solved & \cite{CHALUPKA2024} \\\hline
 $ax^5+by^3=cz^{11}$ & $(a,b,c) \in \cal{C}$, see \eqref{eq:putz} & Solved & \cite{ThesisPutz} 
\\\hline
\end{longtable}
\footnotetext{{\url{https://mathoverflow.net/questions/480114}}}

\section*{Acknowledgements}
The authors sincerely thank Frits Beukers, Andrew Bremner, David Zureick-Brown, Tim Browning, Carmen Bruni, Pedro Jos\'e Cazorla Garc\'\i a, Andrzej D\c abrowski, Rainer Dietmann, Jordan Ellenberg, Victor Flynn, Nuno Freitas, Andrew Granville, Alain Kraus, Florian Luca, Lo\"ic Merel, Oana P\u{a}durariu, Kenneth Ribet, Jeremy Rouse, Samir Siksek, Sankar Sitaraman, G\"okhan Soydan, Alain Togb\'e, Lucas Villagra Torcomian, Soroosh Yazdani, Konstantine Zelator and Zhongfeng Zhang for their responses, feedback and corrections on earlier versions of this paper. This research used the ALICE High Performance Computing facility at the University of Leicester. We also thank the referees for their useful comments to improve the quality of this paper.

\bibliographystyle{habbrv}

\bibliography{genfermat.bib}

\enddocument